\pdfoutput=1
\documentclass{article}
\usepackage[T1]{fontenc}
\usepackage[UKenglish]{babel}

\usepackage{Factory, Math, Theorema, Classic, Styl}

\newcommand{\ieme}{$^{\text{e}}$}
\newcommand{\no}{n$^{\text{o}}$}

\makeatletter
    \def\ordinal#1{\@ordinal{\@nameuse{c@#1}}}   
    \def\@ordinal#1{\ifcase#1 Noughth\or First\or Second\or Third\or Fourth\else\@iordinal{#1}\fi}
    \def\@iordinal#1{\ifcase#1\or \or \or \or \or Fifth\or f\or g\or h\or i\or j\or%
       k\or l\or m\or n\or o\or p\or q\or r\or s\or t\or u\or v\or w\or x\or%
       y \or z\else\@ctrerr\fi}
\makeatother

\renewcommand{\thesection}{\arabic{section}} 
\renewcommand{\thesubsection}{\S\arabic{subsection}} 
\titleformat{\section}[display]{\centering\scshape\normalsize}{\small Section the \ordinal{section}.}{0pt}{} 
\titleformat{\subsection}[runin]{\bfseries}{\thesubsection.}{1ex}{}[.]

\titlecontents{section}[1.8em]{\addvspace{1pc}\bfseries}
{\contentslabel[\thecontentslabel.\hspace{1ex}]{1.8em}}
{\hspace*{-1.8em}}
{\titlerule*[1pc]{}\contentspage}

\titlecontents{subsection}
[4.1em]
{}
{\contentslabel[\thecontentslabel.\enspace\hfill]{2.3em}}
{}
{\titlerule*[1pc]{.}\contentspage}

\usepackage[all, 2cell]{xy} \UseAllTwocells 
\newcommand{\upar}{ \ar@<+.4ex> }			\newcommand{\downar}{ \ar@<-.4ex> }
\newcommand{\curvear}{\ar@/^1pc/} 			

\newcommand{\B}{\mathbf{B}}

\newcommand{\pass}[3]{ \xymatrix{#1 \ar[r]^{#3} & #2}} 
\newcommand{\spass}[3]{ \left\{ \xymatrix{#1 \ar[r]^{#3} & #2} \right\} } 

\newcommand{\ssap}[3]{ \xymatrix{#1 & \ar[l]_{#3} #2}} 
\newcommand{\bssap}[3]{ \begin{bmatrix} \xymatrix{#1 & \ar[l]_{#3} #2} \end{bmatrix}} 

\newcommand{\pcom}[3]{#1(#2^{[#3]})} 
\newcommand{\col}[2]{\begin{bmatrix} #1 \\ #2 \end{bmatrix}}
\newcommand{\smallcol}[2]{\left[ \begin{smallmatrix} #1 \\ #2 \end{smallmatrix} \right]}

\newcommand{\Mod}{\mathfrak{Mod}} 		\newcommand{\XMod}{\mathfrak{XMod}}	\newcommand{\BMod}{\Mod_\B}
\newcommand{\Sets}{\mathfrak{Set}}		\newcommand{\MSet}{\mathfrak{MSet}} 
\newcommand{\Laby}{\mathfrak{Laby}}		\newcommand{\Sur}{\mathfrak{Sur}}		\newcommand{\NQO}{\widehat{\Sur}}
\newcommand{\MFun}{\Fun(\XMod,\Mod)}		 \newcommand{\HPol}{\mathfrak{Hom}}
\newcommand{\Num}{\mathfrak{Num}} 		\newcommand{\QPol}{\mathfrak{QHom}}
\newcommand{\CAlg}{\mathfrak{CAlg}}

\newcommand{\de}{\diamond} 			\DeclareMathOperator*{\De}{\lozenge} 

\DeclareMathOperator{\Lin}{\mathrm{Lin}}
\DeclareMathOperator{\ce}{\mathrm{cr}}	

\usepackage{fancyhdr}						
\begin{document}

\lefthyphenmin=2 \righthyphenmin=2

\title{\Large\bfseries THE COMBINATORICS OF \\ POLYNOMIAL FUNCTORS  \\[2.2ex] \rule{10em}{.7pt} \\[1ex] }
\author{\normalsize\sc Qimh Richey Xantcha%
\thanks{\textsc{Qimh Richey Xantcha}, Uppsala University: \texttt{qimh@math.uu.se}}}
\date{\normalsize\textit{\today}}
\maketitle

\bigskip 
\epigraph{\begin{vverse}[Tusen vägar med irrande svek: Knappt mäktar han sjelf att]
Dedalus genom sin konst och sitt snille vida beryktad \\
Bygde det opp; han förvirrar de ledande märken och ögat \\
I villfarelse för ibland skiljaktiga vägar. \\
Så på de Frygiska fält, man ser den klara Meandros \\
Leka. I tveksamt lopp han rinner och rinner tillbaka, \\
Möter sig ofta sjelf och skådar sin kommande bölja, \\
Och nu till källan vänd, nu åt obegränsade hafvet, \\
Rådvill öfvar sin våg. Så fyllas af Dedalus äfven \\
Tusen vägar med irrande svek: Knappt mäktar han sjelf att \\
Hitta till tröskeln igen. Så bedräglig han boningen danat. \par
\vattr{Ovidius, \emph{Metamorphoses}}
\end{vverse}}

\bigskip 
\renewcommand{\abstractname}{\mdseries\scshape Argument}
\begin{abstract} \noindent\footnotesize
We propose a new description of Endofunctors of Module Categories, based upon a 
combinatorial category comprising finite sets and so-called \emph{mazes}. Polynomial and numerical 
functors both find a natural interpretation in this frame-work. 

Since strict polynomial functors, according to the work of Salomonsson, are encoded by multi-sets, the 
two strains of functors 
may be compared and contrasted through juxtaposing the respective combinatorial structures, leading to the 
\emph{Polynomial Functor Theorem}, giving an effective criterion for when a numerical (polynomial) functor 
is strict polynomial.  

\begin{description}
\item[\mdseries\textsc{2010 Mathematics Subject Classification.}]  Primary: 16D90. Secondary: 13C60, 18A25.
\end{description}
\end{abstract}

\bigskip

\noindent Polynomial and strict polynomial functors abound. Recent years have seen scholars 
vying for  
adequate descriptions of their different flavours, 
with their main focus on quadratic ones, oft-times being the only manageable case. 
Let us, for the moment, be contented with mentioning Baues and Pirashvili's classification \cite{BP} 
of quadratic functors on Groups, and Hartl and Vespa's account \cite{HV} 
of quadratic functors from Pointed Categories to Abelian Groups. 
Also, it will not be superfluous to note Hartl, Pirashvili, and Vespa's recently completed classification 
\cite{HPV} of polynomial functors from Free Algebras 
over set operads (which includes the cases of abelian groups and groups) to Abelian Groups.    

We shall be concerned with classifying endofunctors on Module Categories, with 
a particular eye on polynomial ones. 
Work in this direction was presumably initiated by Baues, with his classification \cite{Baues} of quadratic functors, 
and then widely expanded upon by Baues, Dreckman, Franjou, and Pirashvili in their now classical tract \cite{BDFP}, 
in which polynomial endofunctors of Abelian Groups are given an accurate characterisation. 
To this end, the authors deploy the category of finite sets and surjections, an ingenious feat, 
but not without certain draw-backs. Firstly, their scheme fails to bring about a classification of 
\emph{all} abelian group functors, polynomial or not. 
This is because, at the core of their argument, lies a structure theorem 
for functors on Commutative \emph{Monoids}. 
Passing to polynomial functors, it will be recalled, 
has the curious effect of eradicating the distinction between monoids and groups.  
Secondly, attempts to 
generalise their argument to arbitrary module categories will encounter quandaries, as it is not immediately 
clear what category one should employ in lieu of Surjections. 

We propose to resolve both of these difficulties by introducing the \emph{Labyrinth Category} $\Laby$, consisting of 
finite sets with \emph{mazes} as the arrows between them. A \emph{maze} from the set $\{a,b,c,d\}$ to the set $\{w,x,y,z\}$ (say) is something like the following diagram:
{\small $$ \begin{bmatrix} \xymatrixcolsep{2.5em} \xymatrixrowsep{1.5em} 
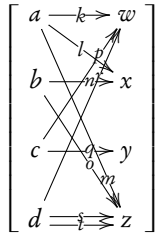
\xymatrix{
a \ar[r]|{k} \ar[dr]|{l} \ar[dddr]|(0.8){m} 	& w \\
b \ar[r]|(0.6){n} \ar[ddr]|(0.6){o}	    	& x \\
c \ar[uur]|(0.7){p} \ar[r]|(0.6){q} 	    	& y \\
d \ar[uuur]|(0.7){r} \upar[r]|{s} \downar[r]|{t}& z  
} \end{bmatrix} $$}
where $k,l,m,n,o,p,q,r,s,t$ are scalars. 
The reason for applying the name \emph{maze} to such a contraption should not be difficult to divine. 

The Labyrinth Category successfully captures the skeletal structure of the 
whole module category and enables the reduction of a functor to  
a significantly smaller collection of data. Polynomial functors find a natural interpretation in this frame-work, 
as do numerical functors. (The honourable reader is referred to the preliminary section for definitions.) 
Our main theorems along this line are as follows. The base ring is betokened by $\B$ and its category of right modules by 
 $\Mod=\BMod$. The symbol $\XMod$ denotes the subcategory of finitely generated, 
free modules. Functors from $\Laby$ are always assumed linear, while functors on the module category may 
be of quite an arbitrary nature.

\begin{inttheorem}[\ref{T: LoF}: Labyrinth of Fun]
The functor $$ \Phi\colon \MFun\to \Lin(\Laby,\Mod), $$
where $\Phi(F)\colon \Laby\to\Mod$ maps the set $X$ to the cross-effect of rank $X$ of $F$, evaluated on $\B$;
is an equivalence of categories. 
\end{inttheorem}

\begin{inttheorem}[\ref{T: LoF Pol}]
The module functor $F$ is polynomial of degree $n$ if and only if 
$\Phi(F)$ vanishes on sets with more than $n$ elements.
\end{inttheorem}

\begin{inttheorem}[\ref{T: LoF Num}] Assuming a binomial base ring, 
the module functor $F$ is numerical of degree $n$ if and only if it suffices to specify the action 
of $\Phi(F)$ on \textbf{pure} mazes; i.~e., mazes carrying only the label $1$. 
\end{inttheorem}

Roughly speaking, mazes should be thought of as deviations, and the latter theorem substantiates the known fact 
that, for numerical 
functors, co\hyp efficients from inside the argument may be brought out front, 
as testified by Theorem 10 of \cite{PF}.
The quotient category encoding numerical functors will be denoted by $\Laby_n$, so that: 
\emph{Numerical functors, of degree $n$, are equivalent to linear functors from $\Laby_n$ to Modules.}

We interpose here, referring to Theorem \ref{T: BDFP} below for the precise statement, 
that the labyrinth description will indeed be equivalent to Baues, Dreckmann, Franjou, and Pirashvili's 
approach \cite{BDFP} using surjections, 
contingent on the hypotheses that: (1) polynomiality have been duly assumed, and (2) the base ring be $\Z$. 

The further question, as to whether strict polynomial functors admit a similar formatting, 
found its complete solution in  Salomonsson's doctoral thesis \cite{Pelle}. 
Denoting the category of multi-sets of cardinality $n$ by $\MSet_n$, 
the main result that shall concern us  
may, slightly paraphrased, be summarised thus: 
\emph{Homogeneous functors, of degree $n$, are equivalent to linear functors from $\MSet_n$ to Modules.} 
(See Theorem \ref{T: LoF Hom}.) 

According to Theorem \ref{T: BZ} below, there is an isomorphism 
$$ {}_\B\Laby_n \cong \B\otimes_\Z {}_\Z\Laby_n, $$
whose heuristical significance it may not be superfluous to point out.
Namely, the category of numerical functors over an \emph{arbitrary} binomial ring will be found 
``identical in structure'' to the category of numerical --- or polynomial --- 
functors over $\Z$, for
\begin{align*}
{}_\B\Num_n 
&\cong \Lin_{\B}({}_\B\Laby_n,\Mod) \cong \Lin_{\B}(\B\otimes_\Z {}_\Z\Laby_n,\Mod) \\
&\cong \B\otimes_\Z  \Lin_{\Z}({}_\Z\Laby_n,\Mod) \cong \B\otimes_\Z {}_\Z\Num_n.
\end{align*}
The corresponding result for homogeneous functors holds more trivially, for here it is an inherence of the 
construction that 
$$ {}_\B\MSet_n \cong \B\otimes_\Z {}_\Z\MSet_n. $$

One obvious mode of relating the two strains of functors --- numerical and strict polynomial --- 
will be to juxtapose their combinatorics. 
To this end, we establish the existence of a functor 
$$ A_n\colon \Laby_n\to\MSet_n, $$ called the \emph{Ariadne functor}, by the aid of which two principal 
results are procured: 

\begin{inttheorem}[\ref{T: Ariadne}] Pre-composition with the Ariadne functor begets 
the forgetful functor from Homogeneous to Numerical Functors. \end{inttheorem}

\begin{inttheorem}[\ref{T: PFT}: The Polynomial Functor Theorem] Let $F$ be a numerical functor of degree $n$, corresponding 
to the labyrinth module $H\colon\Laby_n\to\Mod$. 
Then $F$ may be given a homogeneous structure of degree $n$ if and only if the following criteria are met:
\bnum
\item $F$ is \textbf{quasi-homogeneous}, that is, it satisfies the equation $$ F(r\alpha) = r^nF(\alpha), $$
for any $r\in\Q\otimes_\Z \B$ and homomorphism $\alpha$ (see Definition \ref{D: QH}).
\item $H$ admits a factorisation through $\Laby^{\oplus n}$ (see Definition \ref{D: Laby+}).
\enum \hfill 
\end{inttheorem}

An equipollent variant is Theorem 23 of \cite{PF}, stated using the language of modules. 

The construction of $\Laby^{\oplus n}$ should be of considerable interest, as it supplies the key to appreciating 
the exact obstructions for polynomial functors to be strict polynomial.  

This research was carried out while a graduate student at Stockholm University under the eminent supervision 
of Prof.~Torsten Ekedahl.

\tableofcontents

\setcounter{section}{-1} \section{Polynomial Functors}

We begin by reviewing the basic theory of Polynomial and Strict Polynomial Functors. 

For the entirety of this article, $\B$ shall denote a fixed base ring of scalars.
All modules, homomorphisms, and tensor products shall be taken over this $\B$, unless other-wise stated. 
We let $\Mod={}_\B \Mod$ betoken the category of (unital) modules over this ring, and we denote
by $\XMod$ be the category\footnote{The 
letter \emph{X} herein is intended to 
suggest ``eXtra nice modules''.} of those modules that are 
finitely generated and free. 

When $A$ and $B$ are linear categories (enriched over $\Mod$), 
the symbols $\Fun(A,B)$ and $\Lin(A,B)$ denote the corresponding
categories of functors and linear functors, respectively.

We shall also use the standard notation $ [n] = \{1,\dots,n\}$. 

\subsection{Polynomial Functors} A \textbf{module functor} is a functor $ \XMod \to \Mod$ --- 
and, so as to avoid any misunderstandings, we duly emphasise that linearity will \emph{not} be assumed. 
We shall be wholly content to consider such restricted functors exclusively, for   
a functor defined on the subcategory $\XMod$ always 
has a canonical well-behaved extension 
to the whole module category $\Mod$. The details are spelled out in \cite{PF}. 

Let us recall the classical notions of polynomiality. 
The subsequent definitions made their first appearance in print, albeit somewhat implicitly, 
in Eilenberg \& Mac Lane's monumental article \cite{EM}, sections 8 and 9: 

\bdf 
Let $\phi\colon M\to N$ be a map of modules. The $n$'th 
\textbf{deviation} of $\phi$ is the map 
$$ \phi(x_1\de\cdots\de x_{n+1}) = \sum_{I\subseteq[n+1]} (-1)^{n+1-|I|} \phi\left(\sum_{i\in I} x_i\right) $$
of $n+1$ variables. \edf

\bdf
The map $\phi\colon M\to N$ is \textbf{polynomial} of degree $n$ if its $n$'th deviation vanishes: 
$$\phi(x_1\de\cdots\de x_{n+1}) =0 $$
for any $x_1,\dots,x_{n+1}\in M$. \edf

\bdf The functor $F\colon \XMod \to \Mod$ is \textbf{polynomial} of degree (at most) $n$ 
if every arrow map $$F\colon \Hom(M,N) \to \Hom(F(M),F(N))$$ is. \edf 

Classically, module functors have been analysed in terms of their \emph{cross-effects}. 

\bdf Let $F$ be a module functor and let $k$ be a natural number. 
The \textbf{cross-effect} of rank $k$ is the multi-functor 
$$ \ce_k F(M_1,\dots,M_k) = \Im F\left(\pi_1\de\cdots\de\pi_k \right), $$
where $\pi_i\colon \bigoplus M_i \to \bigoplus M_i$ denote the canonical projections. \edf

\begin{theorem*}[The Cross-Effect Decomposition (\cite{EM}, Theorem 9.1)] 
$$F(M_1\oplus\cdots \oplus M_k)=\bigoplus_{I\subseteq [k]} \ce_I F\left( (M_i)_{i\in I}\right). $$ \end{theorem*}

A functor is polynomial of degree $n$ if and only if the cross-effects of rank exceeding $n$ vanish.

\subsection{Numerical Functors}
Next, let us suppose the base ring $\B$ to be \emph{binomial}\footnote{This is a \emph{numerical ring} in the 
terminology \cite{TE} of Ekedahl. 
For the equivalence of the two notions, a proof is offered in \cite{BR}.} in the sense of Hall (\cite{Hall}); 
that is, commutative, unital, and in the possession of binomial co-efficients.   
Examples include the ring of integers, as well as all $\Q$-algebras. 

\begin{refdefinition}[\cite{PM}, Definition 5] The map $\phi\colon M\to N$ is \textbf{numerical} of degree (at most) $n$ if it satisfies the following two equations: 
\begin{gather*}
\phi(x_1\de\cdots\de x_{n+1}) =0, \qquad x_1,\dots,x_{n+1}\in M; \\
\phi(rx) = \sum_{k=0}^n \binom{r}{k}\phi\left(\De_k x\right), \qquad r\in\B,\ x\in M.
\end{gather*} \end{refdefinition}

\begin{refdefinition}[\cite{PF}, Definition 10] 
The functor $F\colon \XMod \to \Mod$ is  \textbf{numerical} of degree (at most) $n$ if every arrow map 
$$F\colon \Hom(M,N) \to \Hom(F(M),F(N))$$ is. \end{refdefinition}

We denote the abelian category of numerical functors of degree $n$ by $\Num_n$. 

\begin{refdefinition}[\cite{PF}, Definition 11] \label{D: QH} 
The numerical functor $F$ is \textbf{quasi\hyp homogeneous} of degree $n$ if the extension 
functor\footnote{It will be recalled that binomial rings are torsion-free.} 
$$ F \colon \Q\otimes_\Z\XMod \to \Q\otimes_\Z \Mod $$
satisfies the equation 
$$ F(r\alpha) = r^nF(\alpha), $$
for any $r\in\Q\otimes_\Z \B$ and homomorphism $\alpha$. \end{refdefinition}

The category of quasi-homogeneous functors of degree $n$ will be denoted by the symbol $ \QPol_n$.
Being quasi-homogeneous is an obvious necessary condition for a functor to 
admit a homogeneous structure (see below). The Polynomial Functor Theorem below will provide a sufficient condition. 

\subsection{The Deviations} Since the Labyrinth Category formalises deviations, it will be deemed pivotal, for a proper 
understanding of the theory to follow, to acquire some cognisance of the machinery of the latter. 
To this end, we establish the \emph{Deviation Formula}, which seems to have made its first appearance in \cite{X}. 

For (multi-)sets $A$ and $B$, we shall write $K\sqsubseteq A\times B$ to indicate that $K\subseteq A\times B$ and 
that both canonical projections are surjective.

\blem Let $m$ and $n$ be natural numbers, and let $L\subseteq [m]\times [n]$. Then 
$$ \sum_{L \subseteq K \sqsubseteq [m]\times [n]} (-1)^{|K|} = 0, $$
unless $L$ is of the form $P\times Q$, for some $P\subseteq [m]$, $Q\subseteq [n]$.\elem

\bpr If $L$ is not of the given form, there exists an $(a,b)$ 
which is not in $L$, but such that some $(a,j)$ and some $(i,b)$ are in $L$. Then, for any 
$K\subseteq [m]\times [n]$ containing $(a,b)$, the set $K$ will satisfy the given inclusions 
if and only if $K\setminus \{(a,b)\}$ does. 
Because the cardinalities of these sets differ by $1$, the corresponding terms in the 
sum will have opposite signs, and hence cancel. \epr

\blem Let $m$, $n$, $p$, and $q$ be natural numbers. Then $$\sum_{[p]\times [q] \subseteq K \sqsubseteq [m]\times [n]} (-1)^{|K|} = (-1)^{m+n+p+q+pq}.$$ \elem

\bpr Let $Y(m,n,k)$ denote the number of sets $K$ of cardinality $k$ satisfying 
$$[p]\times [q] \subseteq K \sqsubseteq [m]\times [n].$$
The formula is evidently true for $m=p$ and $n=q$, for then $Y(p,q,pq)=1$, and all other $Y(p,q,k)=0$. 

We now proceed by recursion. 
Consider the pair $(m,n)\in [m]\times [n]$. The sets $K$ containing $(m,n)$ will fall into two classes: 
those where $(m,n)$ is \emph{mandatory} in order to satisfy $K \sqsubseteq [m]\times [n]$, and those where it is not. 
For a $K$ in the latter class, removing $(m,n)$ will yield another set counted in the sum above, but of cardinality decreased by $1$. Since these two types of sets 
exactly pair off, with opposing signs, their contribution to the given sum is $0$. 

Consider then those $K$ of which $(m,n)$ \emph{is} a mandatory element. They fall into three categories:
\begin{itemize}
 \item Some $(m,j)\in K$, for $1\leq j\leq n-1$, but no $(i,n)\in K$, for $1\leq i\leq m-1$. The number of such sets is $Y(m,n-1,k-1)$. 
 \item No $(m,j)\in K$, for $1\leq j\leq n-1$, but some $(i,n)\in K$, for $1\leq i\leq m-1$. The number of such sets is $Y(m-1,n,k-1)$. 
 \item No $(m,j)\in K$, for $1\leq j\leq n-1$, and no $(i,n)\in K$, for $1\leq i\leq m-1$. The number of such sets is $Y(m-1,n-1,k-1)$. 
\end{itemize}
Induction yields
\begin{align*}
& \sum_k (-1)^k Y(m,n,k) \\ &= \sum_k (-1)^k \big( Y(m,n-1,k-1) + Y(m-1,n,k-1) + Y(m-1,n-1,k-1) \big) \\
&= -\big( (-1)^{m+n-1+p+q+pq} + (-1)^{m-1+n+p+q+pq} +(-1)^{m-1+n-1+p+q+pq} \big) \\
& = (-1)^{m+n+p+q+pq}, 
\end{align*}
as desired.  \epr

\begin{theorem}[The Deviation Formula] \label{T: Deviation Formula} 
For a module functor $F$ and homomorphisms $\alpha_1,\dots,\alpha_m$, $\beta_1,\dots,\beta_n$, the following equation holds:
$$ F(\alpha_1\de\cdots\de\alpha_m) \circ F(\beta_1\de\cdots\de\beta_n) = 
\sum_{K\sqsubseteq [m]\times [n]} F\left(\De_{(i,j)\in K} \alpha_i\beta_j\right). $$ \end{theorem}

\bpr We have
\begin{align*}
& \sum_{K\sqsubseteq [m]\times [n]} F\left(\De_{(i,j)\in K} \alpha_i\beta_j\right) = 
\sum_{K\sqsubseteq [m]\times [n]} \sum_{L\subseteq K} (-1)^{|K|-|L|} F\left( \sum_{(i,j)\in L} \alpha_i\beta_j \right) \\
&\qquad = \sum_{L\subseteq [m]\times [n]} \sum_{L\subseteq K\sqsubseteq [m]\times [n] } (-1)^{|K|-|L|}  F\left( \sum_{(i,j)\in L} \alpha_i\beta_j \right) \\
&\qquad = \sum_{L\subseteq [m]\times [n]} (-1)^{|L|} F\left( \sum_{(i,j)\in L} \alpha_i\beta_j \right) \sum_{L\subseteq K\sqsubseteq [m]\times [n] } (-1)^{|K|}   \\
&\qquad = \sum_{P\times Q\subseteq [m]\times [n]} (-1)^{|P||Q|} F\left( \sum_{(i,j)\in P\times Q} \alpha_i\beta_j \right) (-1)^{m+n+|P|+|Q|+|P||Q|}   \\
&\qquad = \sum_{P\subseteq [m]} (-1)^{m-|P|} F\left( \sum_{i\in P} \alpha_i \right) \sum_{Q\subseteq [n]} (-1)^{n-|Q|} F\left( \sum_{j\in Q} \beta_j \right)   \\
&\qquad = F(\alpha_1\de\cdots\de\alpha_m) F(\beta_1\de\cdots\de\beta_n).
\end{align*}
In the fifth step the lemmata were used to evaluate the inner sum. \epr

\subsection{Strict Polynomial Functors}
Let us now recall the strict polynomial maps (``lois polynomes'') from the work of Roby and the strict polynomial 
functors introduced by Friedlander and Suslin. The base ring $\B$ may here be taken commutative and unital only. 

\begin{refdefinition}[\cite{Roby}, section 1.2] A \textbf{strict polynomial map} is a natural transformation 
$$\phi\colon M\otimes - \to N\otimes -$$ betwixt functors $\CAlg\to \Sets$, 
where $\CAlg={}_\B\CAlg$ designates the category of commutative, unital algebras 
over the ring $\B$, and $\Sets$ denotes the category of sets. 
\end{refdefinition}

\begin{refdefinition}[\cite{FS}, Definition 2.1] The functor $F\colon \XMod \to \Mod$ is 
 \textbf{strict polynomial} of degree $n$ if the arrow maps
$$ F\colon \Hom(M,N) \to \Hom(F(M),F(N)) $$
have been given a (multiplicative) strict polynomial structure. \label{D: Strict Polynomial} 
The functor is called \textbf{homogeneous} if all its arrow maps are. \end{refdefinition}

Strict polynomial maps and functors decompose as the direct sum of their homogeneous components. 
We shall denote by $\HPol_n$ the abelian category of homogeneous functors of degree $n$.

For strict polynomial, or even strict analytic (see \cite{PF}, \S0), functors, 
the cross-effects may be decomposed further.
Let $F$ be a homogeneous functor. 
The arrow map $$F\colon \Hom(M,N)\to\Hom(F(M),F(N))$$ will factor through the universal homogeneous map 
$$ \Hom(M,N)\to\Gamma^n\Hom(M,N)$$ 
(see \cite{Roby}, where the full details have been expounded upon),
producing a linear map 
$$ F\colon \Gamma^n\Hom(M,N)\to\Hom(F(M),F(N)).$$
This gives meaning to the symbol $\pcom{F}{\alpha}{A}$, when $\alpha_a$, for $a\in\#A$, are homomorphisms and 
$A$ a multi-set of cardinality $n$. (Multi-sets will be discussed below, in Section \ref{S: Multi-Sets}.)

\bdf Let $F$ be a strict polynomial functor, and let $A$ be a multi-set whose support is $[k]$. 
The \textbf{multi-cross-effect} of rank $A$ is the multi-functor 
$$ \ce_A F(M_1,\dots,M_k) = \Im \pcom{F}{\pi}{A}, $$
where $\pi_i\colon \bigoplus M_i \to \bigoplus M_i$ denote the canonical projections. \edf

In particular, the meaning of the symbol $\ce_A F$, when $A$ is a proper set, is unequivocal.

\begin{theorem}[The Multi-Cross-Effect Decomposition] When $F$ is a homogeneous functor of degree $n$, then
$$ \ce_{[k]}F(M_1,\dots,M_k) = \bigoplus_{\substack{\#A=[k]\\ |A|=n}} \ce_A F(M_1,\dots,M_k), $$
and, consequently,  
$$ F(M_1\oplus\dots\oplus M_k) = \bigoplus_{\substack{\#A\subseteq [k] \\ |A|=n}} \ce_A F\left( (M_a)_{a\in\#A} \right). $$  
\eth

\bpr See, for example, \cite{X}, Theorem 10. \epr

\pagebreak

\section{Multi-Sets} \label{S: Multi-Sets}

\subsection{Multi-Sets}
\bdf A \textbf{multi-set} is a pair $$M=(\# M,\deg_M),$$ where $\#M$ is a set and $$\deg_M\colon \#M\to \Z^+$$ is a function, called the \textbf{degree} or \textbf{multiplicity}. 
The underlying set $\#M$ is called the \textbf{support} of $M$. \edf

Obviously, the support is uniquely determined by the degree function, whence it suffices to specify the latter. 

The degree $\deg_M a$ counts the ``number of times $a\in\#M$ occurs in $M$''. 
The degree of the whole multi-set $M$ we define to be 
$$ \deg M = \prod_{x\in\#M} (\deg x)!. $$

\bdf The \textbf{cardinality} of $M$ is $$ |M| = \sum_{x\in M} 1 = \sum_{x\in \#M} \deg x.$$ \edf

The cardinality counts the number of elements \emph{with} multiplicity.
We tacitly assume all multi-sets under discussion to be \emph{finite}, as these are the only ones we will ever need.

\bdf The following operations are defined on multi-sets. 
\bnum 
\item The \textbf{union} $A\cup B$ of $A$ and $B$ is 
$$ 
\deg_{A\cup B} = \max(\deg_A, \deg_B). $$ 

\item 
The \textbf{disjoint union}
$A\sqcup B$ of $A$ and $B$ is 
$$ 
\deg_{A\sqcup B} = \deg_A + \deg_B. $$ 

\item The \textbf{intersection} $A\cap B$ of $A$ and $B$ is 
$$ \deg_{A\cap B} = \min(\deg_A, \deg_B).$$  

\item The \textbf{relative complement} $A\setminus B$ of $B$ in $A$ is 
$$ \deg_{A\setminus B} = \max(\deg_A-\deg_B,0) .$$ 

\item The \textbf{direct product} $A\times B$ of $A$ and $B$ is 
$$ 
\deg_{A\times B} = \deg_A\cdot \deg_B\colon \#A \times \#B \to \Z^+.$$ 
\enum \hfill
\edf

\bdf $A$ is a \textbf{sub-multi-set\footnote{Some scholars would no doubt say 
\emph{multi-subset}.}} of $B$, written $A\subseteq B$, if 
$$\deg_A \leq\deg_B$$ (element-wise inequality). \edf

\subsection{Multations}

\bdf Let $A$ and $B$ be multi-sets of equal cardinality. 
A \textbf{multation} $A\to B$ is
a \emph{sub-multi-set} of $A\times B$ whose multi-set of first co-ordinates equals $A$
and whose multi-set of second co-ordinates equals $B$. \edf 

Informally, a multation pairs off the elements of one multi-set with those of another. 
The degree $\deg_\mu(a,b)$ counts the number of times $a\in A$ is paired off with $b\in B$. 
A multation $A\to B$  may be written as a two-row matrix, with the 
elements of $A$ on top of those of $B$, the way permutations are usually written (indeed, multations 
should be thought of as generalised such). 

Given a multation 
$$ \begin{bmatrix}  a_1 & a_1 & \dots & a_2 & a_2 & \dots \\ b_1 & b_1 & \dots & b_2 & b_2 & \dots \end{bmatrix}, $$
with $m_j$ appearances of the column $\smallcol{a_j}{b_j}$, we shall adopt the perspective of viewing it as a formal product 
$$ \col{a_1}{b_1}^{[m_1]} \col{a_2}{b_2}^{[m_2]} \dots $$
of \emph{divided powers}\footnote{The divided power $z^{[n]}$ should be thought of as $\frac{z^n}{n!}$.}. 

\bex There exist two multations from the multi-set $\{a,a,b\}$ to itself, namely: 
$$
\begin{bmatrix}  a & a & b \\ a & a & b \end{bmatrix} = \col{a}{a}^{[2]}\col{b}{b} \qquad 
\begin{bmatrix}  a & a & b \\ a & b & a \end{bmatrix} = \col{a}{a} \col{a}{b} \col{b}{a}.
$$
\eex

\subsection{The Multi-Set Category}

Let $\mu\colon B\to C$ and $\nu\colon A\to B$ be two multations, where $|A|=|B|=|C|=n$. 
Their \textbf{composition} or \textbf{product} $\mu\circ\nu$ is found 
by identifying the co-efficient of $x^{\mu}y^{\nu}$ in the equation 
$$ \left(\sum_{\substack{b\in\#B \\ c\in\#C}} x_{cb}\col{b}{c} \right)^{[n]} \circ \left(\sum_{\substack{a\in\#A \\ b\in\#B}} y_{ba}\col{a}{b} \right)^{[n]} =
\left(\sum_{\substack{a\in\#A \\ b\in\#B \\ c\in\#C}} x_{cb}y_{ba}\col{a}{c} \right)^{[n]}, $$
thought of as an equality of divided powers in the formal variables $x$ and $y$. 

The composition may also be viewed thus. The formal multiplication of columns, given by  
$$\col{b'}{c} \circ \col{a}{b} = \bca \col{a}{c} & \text{if $b=b'$,} \\ \ \ 0 & \text{if $b\neq b'$,} \eca$$
 makes the free module of columns into an algebra $A$. Multation composition is then simply the natural multiplication 
$$ u^{[n]} \circ v^{[n]} = (u\circ v)^{[n]} $$
on the module $\Gamma^n(A)$ of divided $n$'th powers.

\bex To calculate the composition 
$$ \begin{bmatrix}  c & d & d \\ e & e & f \end{bmatrix} \circ \begin{bmatrix}  a & a & b \\ c & d & d \end{bmatrix}, $$
we deploy the equation 
\begin{align*}&\left(x_{ec}\col{c}{e} + x_{ed}\col{d}{e} + x_{fd}\col{d}{f}\right)^{[3]} \circ \left(y_{ca}\col{a}{c} + y_{da}\col{a}{d} + y_{db}\col{b}{d} \right)^{[3]}\\ 
&= \left(x_{ec}y_{ca}\col{a}{e} + x_{ed}y_{da}\col{a}{e} + x_{ed}y_{db}\col{b}{e} + x_{fd}y_{da}\col{a}{f} + x_{fd}y_{db}\col{b}{f} \right)^{[3]}.\end{align*}
Identification of the co-efficient of $x_{ec}x_{ed}x_{fd}y_{ca}y_{da}y_{db}$ yields
$$\begin{bmatrix}  c & d & d \\ e & e & f \end{bmatrix} \circ \begin{bmatrix}  a & a & b \\ c & d & d \end{bmatrix} 
= 2\begin{bmatrix}  a & a & b \\ e & e & f \end{bmatrix} + \begin{bmatrix}  a & a & b \\ e & f & e \end{bmatrix}. $$
                                                      
Similarly, by picking the co-efficients of $x_{ec}x_{ed}^2y_{ca}y_{da}^2$, we shall find 
$$\begin{bmatrix}  c & d & d \\ e & e & e \end{bmatrix} \circ \begin{bmatrix}  a & a & a \\ c & d & d \end{bmatrix} = 3\begin{bmatrix}  a & a & a \\ e & e & e \end{bmatrix}. $$
\eex

There is a simple, combinatorial rule for calculating the composition. 
Namely, the composition of two \emph{ordinary} products (\emph{not} divided powers) of columns 
is found by ``summing over all possibilities of composing them'':
$$\left(\col{b_1}{c_1} \cdots \col{b_n}{c_n}\right) \circ \left(\col{a_1}{b_1} \cdots \col{a_n}{b_n} \right)= 
\sum_{\sigma} \left(\col{a_1}{c_{\sigma(1)}} \cdots \col{a_n}{c_{\sigma(n)}} \right) , $$
where the sum is to be taken over all permutations $\sigma\colon [n]\to[n]$ such that $b_j=b_{\sigma(j)}$ for all $j$. 
We leave it to the reader to check the accuracy of this rule. 

\bex 
Computing according to this device, we find
\begin{align*}
& \begin{bmatrix}  c & d & d \\ e & e & f \end{bmatrix} \circ \begin{bmatrix}  a & a & b \\ c & d & d \end{bmatrix}  = 
\col{c}{e}\col{d}{e}\col{d}{f} \circ \col{a}{c}\col{a}{d}\col{b}{d} \\
&\qquad = \col{a}{e}\col{a}{e}\col{b}{f} + \col{a}{e}\col{a}{f}\col{b}{e} 
= 2\begin{bmatrix}  a & a & b \\ e & e & f \end{bmatrix} + \begin{bmatrix}  a & a & b \\ e & f & e \end{bmatrix}, 
\end{align*}
and similarly 
\begin{align*}
\begin{bmatrix}  c & d & d \\ e & e & e \end{bmatrix} \circ \begin{bmatrix}  a & a & a \\ c & d & d \end{bmatrix} 
&= \frac{1}{2}\col{c}{e}\col{d}{e}\col{d}{e} \circ \frac{1}{2}\col{a}{c}\col{a}{d}\col{a}{d} \\
& = \frac{1}{4} \cdot 2 \col{a}{e}\col{a}{e}\col{a}{e}   
= 3\begin{bmatrix}  a & a & a \\ e & e & e \end{bmatrix}. 
\end{align*}
\eex

The \textbf{identity multation} 
$\iota_A$ of a multi-set $A$ is the multation 
in which every element is paired off with itself. It is clear that composition is associative and that the identity multations act as identities. 

\bdf
The \textbf{$n$'th multi-set category} $\MSet_n$ is a linear category. 
Its objects are formal direct sums of multi-sets of cardinality exactly $n$. 
For two multi-sets $A$ and $B$, the arrow set $\MSet_n(A,B)$ is the free module generated by the multations $A\to B$.
\edf

\subsection{Multi-Sets on Multi-Sets}

Unfortunately, later considerations shall shew the necessity of pushing abstraction up to the next level. 

\bdf Let $M$ be a multi-set. A \textbf{multi-set supported in $M$} is a multi-set supported in the set 
$$ M^\# = \Set{(x,k)\in \#M\times \Z^+ | 1\leq k\leq \deg x}. $$ \edf

When speaking of multi-sets supported in a multi-set $M$, we will let ``degree over $M$'' 
stand for ``degree over $M^\#$''. 

\bex Let $M$ be the set $\{a,b,c\}$ and $N$ the multi-set $\{x,x,y\}$. 
There are three multi-sets with support $M$ and cardinality $4$:
$$ \{a,a,b,c\}, \quad \{a,b,b,c\}, \quad \{a,b,c,c\}. $$
In like wise, since $$N^\# = \{(x,1),(x,2),(y,1)\},$$
there are three multi-sets with support $N$ and cardinality $4$:
$$ \{(x,1),(x,1),(x,2),(y,1)\}, \; \{(x,1),(x,2),(x,2),(y,1)\}, \; \{(x,1),(x,2),(y,1),(y,1)\}.$$
These three multi-sets all have degree $2$ over $N$.

We shall usually, when deemed suitable, identify these multi-sets with the collection
$$ \{x,x,x,y\}, \quad \{x,x,x,y\}, \quad \{x,x,y,y\}.$$
\eex

\section{Mazes}

\subsection{Mazes}
Let $X$ and $Y$ be finite sets. A \textbf{passage} from $x\in X$ to $y\in Y$ is a (formal) arrow 
$p$ from $x$ to $y$, labelled with 
an element of $\B$, denoted by $\overline p$. This we write as 
$$p\colon x\to y \quad \text{ or} \quad  \pass{x}{y}{\overline{p}}.$$

\bdf
A \textbf{maze} from $X$ to $Y$ is a multi-set of passages from $X$ to $Y$. It is required that there be at 
least one passage leading
from every element of $X$ and at least one passage leading to every element of $Y$. 
(We, so to speak, wish to prevent \emph{dead ends} from forming.) 
\edf

Because a maze is a multi-set, there can (and, in general, will) be multiple passages between any two given elements. 

It is perfectly legal to consider the \textbf{empty maze} $\emptyset\to\emptyset$. It is the only maze into or out 
of $\emptyset$, and the only maze possessing no passages. 

\bdf
We say $P\colon X\to Y$ is a \textbf{submaze} of $Q\colon X\to Y$ if $P\subseteq Q$ as multi-sets. This will be denoted 
by the symbol $P\leq Q$.
\edf

\subsection{The Labyrinth Category}

Passages $p\colon y\to z$ and $q\colon x\to y$ are said to be \textbf{composable}, 
seeing that one ends where the other begins. 

\bdf
If $P\colon Y\to Z$ and $Q\colon X\to Y$ are mazes, their \textbf{cartesian product} 
$P\boxtimes Q$ is the multi-set of all pairs of composable passages: 
$$ \xymatrixcolsep{1pc}
P\boxtimes Q = \Set{ \left( \bssap{z}{y}{p} , \bssap{y}{x}{q} \right) | \bssap{z}{y}{p} \in P \och \bssap{y}{x}{q} \in Q } .
$$
\edf

Recall that we, for a sub-multi-set $U\subseteq P\boxtimes Q$, write $U\sqsubseteq P\boxtimes Q$ 
to indicate that the projections on $P$ and $Q$ are both surjective.  
Such a set $U$ can be naturally interpreted as a maze itself, viz.:
$$ \xymatrixcolsep{1pc}
\Set{ \bssap{z}{x}{pq}  | \left( \bssap{z}{y}{p} , \bssap{y}{x}{q} \right)\in U }
$$
(observe the order in which $p$ and $q$ occur).
The surjectivity condition on the projections will prevent dead ends from forming. 
Henceforth, this identification will be made without comment.

\bex
Consider the two mazes 
$$ \xymatrixrowsep{0pc} \xymatrixcolsep{0.5pc}
P=\begin{bmatrix}  \xymatrix{
x  \\
& z \ar[ul]_a \ar[dl]^b  \\
y 
} \end{bmatrix} , \qquad
Q = \begin{bmatrix} \xymatrix{
& x \ar[dl]_c   \\
z   \\
& y \ar[ul]^d 
} \end{bmatrix}. $$
Their cartesian product is 
\begin{align*} 
P\boxtimes Q &= \xymatrixcolsep{1pc} \bigg\{ \left( \bssap{x}{z}{a} , \bssap{z}{x}{c} \right), \left(\bssap{y}{z}{b} , \bssap{z}{x}{c} \right), \\
& \qquad \xymatrixcolsep{1pc} \left( \bssap{x}{z}{a} , \bssap{z}{y}{d} \right), \left(\bssap{y}{z}{b} , \bssap{z}{y}{d} \right) \bigg\},
\end{align*}
which we identify with the maze 
$$
\xymatrixcolsep{2.5em} \xymatrixrowsep{2.5em}
\begin{bmatrix}  \xymatrix{
x & \ar[l]_(0.7){ac} \ar[dl]|(0.7){bc}  x \\
y & \ar[ul]|(0.7){ad} \ar[l]^(0.7){bd} y
} \end{bmatrix}.
$$ 
\eex

We now define the composition of two mazes. As for multi-sets, this composition will not in general be a maze, 
but rather a \emph{sum} of mazes, 
and living in the free module generated by those. 

\bdf The \textbf{composition} or \textbf{product} of the mazes $P$ and $Q$ is the formal sum 
$$ P\circ Q = \sum_{U\sqsubseteq P\boxtimes Q} U . $$ \edf

\bex Let $P$ and $Q$ be as above. Their composition is  
\begin{align*} 
PQ &=
\xymatrixcolsep{2.5em} \xymatrixrowsep{2.5em}
\begin{bmatrix}  \xymatrix{
x & \ar[l]_(0.7){ac} \ar[dl]|(0.7){bc}  x \\
y & \ar[ul]|(0.7){ad} \ar[l]^(0.7){bd} y
} \end{bmatrix}
+\begin{bmatrix}  \xymatrix{
x & \ar[l]_(0.7){ac}  x \\
y & \ar[l]^(0.7){bd} y
} \end{bmatrix}
+\begin{bmatrix}  \xymatrix{
x & \ar[dl]|(0.7){bc}  x \\
y & \ar[ul]|(0.7){ad}  y
} \end{bmatrix} \\
&+\begin{bmatrix}  \xymatrix{
x &  \ar[dl]|(0.7){bc}  x \\
y & \ar[ul]|(0.7){ad} \ar[l]^(0.7){bd} y
} \end{bmatrix}
+\begin{bmatrix}  \xymatrix{
x & \ar[l]_(0.7){ac}   x \\
y & \ar[ul]|(0.7){ad} \ar[l]^(0.7){bd} y
} \end{bmatrix}
+\begin{bmatrix}  \xymatrix{
x & \ar[l]_(0.7){ac} \ar[dl]|(0.7){bc}  x \\
y & \ar[l]^(0.7){bd} y
} \end{bmatrix}
+\begin{bmatrix}  \xymatrix{
x & \ar[l]_(0.7){ac} \ar[dl]|(0.7){bc}  x \\
y & \ar[ul]|(0.7){ad}  y
} \end{bmatrix} .
\end{align*} \eex

That composition is associative follows from the observation that $(PQ)R$ and $P(QR)$ both equal
$$ \sum_{W\sqsubseteq P\boxtimes Q\boxtimes R} W $$
(surjective projections on all three factors). There exist \textbf{identity mazes} 
$$ 
I_X = \bigcup_{x\in X} \spass{x}{x}{1} . $$

\bdf The \textbf{Labyrinth Category} $\Laby$ is a linear category. Its objects are formal direct sums of finite sets. 
Given two sets $X$ and $Y$, the arrow set $\Laby(X,Y)$ is the module generated by the mazes $X\to Y$, with 
the following relations imposed, for any multi-set $P$ of passages: 
\brom
\item $$
\Bigg[ P\cup \spass{\ast}{\ast}{0}  \Bigg] = 0.
$$ 
\item 
\begin{align*}
& \Bigg[ P\cup \spass{\ast}{\ast}{a+b}  \Bigg] = \\
& \Bigg[ P\cup \spass{\ast}{\ast}{a}  \Bigg] + \Bigg[ P\cup \spass{\ast}{\ast}{b}  \Bigg] + 
\Bigg[ P\cup\left\{ \xymatrix{\ast\upar[r]^a \downar[r]_b & \ast } \right\} \Bigg].
\end{align*}
\erom \hfill \edf

The second axiom may be generalised by means of mathematical induction to yield the following elementary formul\ae. 

\bth In the Labyrinth Category, the following two equations hold: \label{T: Maze formulae}
\begin{align*}
\Bigg[ P\cup \spass{\ast}{\ast}{\sum_{i=1}^n a_i}  \Bigg] &= 
\sum_{\emptyset\subset I\subseteq [n]} \Bigg[ P\cup \bigcup_{i\in I} \spass{\ast}{\ast}{a_i} \Bigg] \\
 \Bigg[ P\cup \bigcup_{i=1}^n \spass{\ast}{\ast}{a_i} \Bigg] &=   
\sum_{I\subseteq [n]} (-1)^{n-|I|}  \Bigg[ P\cup \spass{\ast}{\ast}{\sum_{i\in I} a_i} \Bigg].
\end{align*}
\eth

\subsection{The Quotient Labyrinth Categories}

The theory developed unto this point makes sense for an arbitrary base ring. 
In order to construct quotient categories of $\Laby$, we need the 
assumption that $\B$ be binomial.  

When $A$ is a maze, let $I_A$ denote the maze 
$$ I_A = \bigcup_{[p\colon x\to y]\in A} \spass{x}{y}{1}, $$
in which all passages of $A$ have been re-assigned the label $1$. 

\bdf The category $\Laby_n$ is the quotient category obtained from $\Laby$ when imposing the 
following relations, for any maze $P$: 
\brom 
\item[III.] $$\displaystyle P=0, \quad\text{whenever $|P|>n$.} $$
\item[IV.] $$\displaystyle P= \sum_{\#A=P} \prod_{p\in P} \binom{\overline p}{\deg_A p} I_A.$$ 
(The sum is in fact finite, owing to the previous axiom.) \erom \hfill\edf

\bex When $n=3$, an instance of the fourth axiom is the following:
$$ \xymatrixcolsep{.8pc} 
\begin{bmatrix}  \xymatrix{ \ast \ar@/^1pc/[r]^a \ar@/^-1pc/[r]_b  & \ast  }  \end{bmatrix}
\xymatrixcolsep{.8pc} = 
\binom a1 \binom b1  \begin{bmatrix}  \xymatrix{  \ast \ar@/^1pc/[r]^1 \ar@/^-1pc/[r]_1  & \ast }\end{bmatrix}  + 
\binom a2 \binom b1  \begin{bmatrix}  \xymatrix{  \ast \ar@/^1pc/[r]^1 \ar@/^.3pc/[r]_1 \ar@/^-1pc/[r]_1  & \ast }\end{bmatrix}  + 
\binom a1 \binom b2 \begin{bmatrix}  \xymatrix{  \ast \ar@/^1pc/[r]^1 \ar@/^-.3pc/[r]^1 \ar@/^-1pc/[r]_1  & \ast }\end{bmatrix}.$$ 
\eex 

The original Labyrinth Category $\Laby$ encodes arbitrary module functors.  
Axioms III and IV have been appended in order to encode \emph{polynomial} and 
\emph{numerical} functors, respectively.
It may be shewn that, over the integers, Axiom IV is actually implied by Axiom III, comparable to
how numerical and polynomial functors are equivalent in this setting. 

We shall have reason to impose upon the Labyrinth Category yet another axiom. This will be for encoding 
\emph{quasi-homogeneous} functors (of degree $n$). It will be proved presently that the 
category $\Laby_n$ is free over $\B$, and so, 
 in particular, torsion-free, leading to an inclusion of categories
$$ \Laby_n \subseteq \Q\otimes_\Z \Laby_n. $$

When $P$ is a maze and $a$ is a scalar (in a binomial base ring), denote by $a\boxdot P$ the maze obtained from 
$P$ by multiplying the labels of all passages by $a$: 
$$ a\boxdot P = \Set{ \bssap{y}{x}{ap} | \bssap yxp \in P }. $$

\bdf The category $\Laby^n$ is the quotient category obtained from $\Laby_n$ upon the imposition 
of the following axiom, for any maze $P$:
\brom \item[V.] $$a^nP = a\boxdot P, \quad a\in\Q\otimes_\Z\B.$$ \erom
\hfill \edf

\bex
The fifth axiom considers the ideal generated in $\Q\otimes_\Z \Laby_n$, rather than $\Laby_n$. 
The slightly sharper requirement will first make a difference in degree $4$.  
Dividing out by the ideal \emph{generated in  ${}_\Z\Laby_4$} by elements of the form 
$ a^4P - a\boxdot P$ makes it possible to prove
$$ \xymatrixrowsep{0.7pc} \xymatrixcolsep{1.3pc}
2\begin{bmatrix}  \xymatrix{
 & \ast  \\
\ast \ar@/^.6pc/[ur]^1 \ar[ur]|1 \ar@/_.6pc/[ur]_1 \ar[dr]_1  & \\
& \ast 
} \end{bmatrix} + 
2\begin{bmatrix}  \xymatrix{
 & \ast  \\
\ast \upar[ur]^1 \downar[ur]_1 \upar[dr]^1 \downar[dr]_1  & \\
& \ast 
} \end{bmatrix} 
+ 2\begin{bmatrix}  \xymatrix{
 & \ast  \\
\ast \ar[ur]^1 \upar[dr]^1 \downar[dr]_1  & \\
& \ast 
} \end{bmatrix} = 
12\begin{bmatrix}  \xymatrix{
 & \ast  \\
\ast \ar[ur]^1 \ar[dr]^1 & \\
& \ast 
} \end{bmatrix}, $$
whereas we shall be needing a stronger statement. The full force of Axiom V authorises a division   
 by $2$, thus establishing 
$$ \xymatrixrowsep{0.7pc} \xymatrixcolsep{1.3pc}
\begin{bmatrix}  \xymatrix{
 & \ast  \\
\ast \ar@/^.6pc/[ur]^1 \ar[ur]|1 \ar@/_.6pc/[ur]_1 \ar[dr]_1  & \\
& \ast 
} \end{bmatrix} + 
\begin{bmatrix}  \xymatrix{
 & \ast  \\
\ast \upar[ur]^1 \downar[ur]_1 \upar[dr]^1 \downar[dr]_1  & \\
& \ast 
} \end{bmatrix} 
+ \begin{bmatrix}  \xymatrix{
 & \ast  \\
\ast \ar[ur]^1 \upar[dr]^1 \downar[dr]_1  & \\
& \ast 
} \end{bmatrix} = 
6\begin{bmatrix}  \xymatrix{
 & \ast  \\
\ast \ar[ur]^1 \ar[dr]^1 & \\
& \ast 
} \end{bmatrix}. $$
\eex

\section{The Ariadne and Theseus Functors}

We propose to investigate how the Multi-Set and Labyrinth Categories are related. 
A functor in one direction is readily found; viz.~the \emph{Ariadne functor} 
$$ A_n\colon \Laby\to \MSet_n,$$
so called because it leads the way \emph{out of} the labyrinth. In the case of a binomial base ring, 
it will factor:
$$ A_n\colon \Laby\to\Laby_n\to\Laby^n\to\MSet_n.$$

A minor modification of the category $\Laby^n$ will enable us to define a functor in the reverse direction. 
This is the \emph{Theseus functor $T_n$}, going \emph{into} the labyrinth. 
These two functors are inverse to each other. 
The modifications necessary to undertake on $\Laby_n$ in order to define $T_n$, indicate the precise 
obstructions that may prevent a numerical (polynomial) functor from being strict polynomial. 
Confer the Polynomial Functor Theorem, Theorem \ref{T: PFT}, below.

\subsection{The Ariadne Functor}

For the duration of this section, let $n$ be a fixed natural number.

Let $P$ be a maze. Consider the following sum of multations:
$$ A_n(P) 
= \sum_{\substack{\#A=P \\ |A|=n}} \left( \frac{1}{\deg_P A} \prod_{[p\colon x\to y]\in A} \overline{p}\col{x}{y} 
\right) 
= \sum_{\substack{\#A=P \\ |A|=n}} \left(\prod_{[p\colon x\to y]\in P} \overline{p}^{\deg_A p}\col{x}{y}^{[\deg_A p]} 
\right). $$
This will provide a functor from $\Laby$ to $\MSet_n$, as we now set out to prove.
It is clear that $A_n(P)=0$ if a single passage of $P$ be labelled $0$. Now to shew that 
\begin{align} \label{E: Ariadne Axiom II}
& A_n\Bigg( P\cup \spass{u}{v}{a+b}  \Bigg) = \\ 
& \xymatrixcolsep{1pc}  A_n\Bigg( P\cup \spass{u}{v}{a} \Bigg) + A_n\Bigg( P\cup \spass{u}{v}{b} \Bigg) + 
A_n\Bigg( P\cup\left\{ \xymatrix{u\upar[r]^a \downar[r]_b & v } \right\} \Bigg). \nonumber
\end{align}
Denote, when $A$ is a multi-set of passages with support $P$, 
$$ \mu_A = \frac{1}{\deg_P A} \prod_{[p\colon x\to y]\in A} \overline{p}\col{x}{y}.$$
Since
\begin{align*}
\xymatrixcolsep{1pc} A_n\Bigg( P\cup \spass{u}{v}{a+b}  \Bigg) &= \sum_{m=1}^n \sum_{\substack{\#A=P \\ |A|=n-m}} 
\left(  \mu_A  \cdot (a+b)^m\col{u}{v}^{[m]} \right)\\
\xymatrixcolsep{1pc}  A_n\Bigg( P\cup \spass{u}{v}{a} \Bigg) &= \sum_{m=1}^n \sum_{\substack{\#A=P \\ |A|=n-m}} 
\left(   \mu_A \cdot a^m\col{u}{v}^{[m]} \right) \\
\xymatrixcolsep{1pc} A_n\Bigg( P\cup \spass{u}{v}{b} \Bigg) &= \sum_{m=1}^n \sum_{\substack{\#A=P \\ |A|=n-m}} 
\left(   \mu_A \cdot b^m\col{u}{v}^{[m]} \right) \\
\xymatrixcolsep{1pc} A_n\Bigg( P\cup\left\{ \xymatrix{u\upar[r]^a \downar[r]_b & v } \right\} \Bigg) &= 
\sum_{m=1}^n \sum_{\substack{\#A=P \\ |A|=n-m}} \sum_{\substack{i+j=m \\ i,j\geq 1}}
  \left(   \mu_A \cdot a^i b^j \col{u}{v}^{[i]} \col{u}{v}^{[j]} \right),
\end{align*}
the relation \eref{E: Ariadne Axiom II} follows from the equation 
$$ (x+y)^{[m]} = x^{[m]} + y^{[m]} + \sum_{\substack{i+j=m \\ i,j\geq 1}} x^{[i]}y^{[j]},$$
valid in every divided power algebra. 

Hence $A_n$ gives a well-defined map on the mazes of the Labyrinth Category. We now prove that 
it is, in fact, a functor. 

\bth The formul\ae 
\begin{gather*}
A_n(X) = \bigoplus_{\substack{\#A= X \\ |A|=n}} A \\
A_n(P) = \sum_{\substack{\#A=P \\ |A|=n}} 
\left(\prod_{[p\colon x\to y]\in P} \overline{p}^{\deg_A p}\col{x}{y}^{[\deg_A p]} \right)
\end{gather*}
provide a linear functor $$ A_n\colon \Laby\to \MSet_n.$$ \eth

\bpr Let $P\colon Y\to Z$ and $Q\colon X\to Y$ be two mazes. A straight-forward calculation gives at hand:
\begin{align*} 
& A_n(P) \circ A_n(Q) \\
&\quad = \sum_{\substack{\#A=P \\ |A|=n}} 
\left( \prod_{[p\colon y\to z]\in P} \overline{p}^{\deg_A p}\col{y}{z}^{[\deg_A p]} \right) \circ
\sum_{\substack{\#B=Q \\ |B|=n}} 
\left( \prod_{[q\colon x\to y]\in Q} \overline{q}^{\deg_B q}\col{x}{y}^{[\deg_B q]} \right) \\
&\quad = \sum_{\substack{\#C\sqsubseteq P\boxtimes Q \\ |C|=n}} 
\left( \prod_{[r\colon x\to z]\in \#C} \overline{r}^{\deg_C r}\col{x}{z}^{[\deg_C r]} \right) \\
&\quad = \sum_{R\sqsubseteq P\boxtimes Q} \sum_{\substack{\#C=R \\ |C|=n}} 
\left( \prod_{[r\colon x\to z]\in R} \overline{r}^{\deg_C r}\col{x}{z}^{[\deg_C r]} \right) 
= A_n\left(\sum_{R\sqsubseteq P\boxtimes Q}  R \right) = A_n(PQ).
\end{align*}
The only possibly dubious step here is the third, which follows from the equation 
$$ \left(\sum_{[p\colon y\to z]\in P} \overline{p}\col{y}{z} \right)^{[n]} \circ \left(\sum_{[q\colon x\to y]\in Q} \overline{q}\col{x}{y} \right)^{[n]} =
\left(\sum_{\substack{ [p\colon y\to z]\in P \\ [q\colon x\to y]\in Q }} \overline{p}\overline{q}\col{x}{z} \right)^{[n]}, $$
after noting that restricting attention to monomials $\overline{p}^A\overline{q}^B$ with $\#A=P$ 
and $\#B=Q$ in the left-hand side corresponds to 
considering monomials $(\overline{p}\overline{q})^C$ satisfying the three relations
$$\#C\sqsubseteq P\boxtimes Q, \qquad C_1=A, \qquad C_2=B $$ 
(canonical projections) in the right-hand side. \epr 

\bdf The functor $A_n$ is called the $n$'th \textbf{Ariadne functor}. \edf

\blem When $x$ is an element of a divided power algebra over a binomial ring $\B$, then 
$$ a^m x^{[m]} = \sum_{k=1}^\infty \binom{a}{k} \sum_{\substack{g_1+\cdots+g_k=m \\ g_i\in\Z^+}} x^{[g_1]}\cdots x^{[g_k]}, 
\qquad a\in\B,\ m\in\Z^+. $$
 \elem

\bpr The lemma will be an easy consequence of the formula 
$$ a^m = \sum_{k=1}^\infty \binom ak \sum_{g_1+\cdots+g_k=m} \binom{m}{g_1,\dots,g_k}, $$
established by a simple combinatorial argument.
When $a\in\N$, both sides count the number of ways to colour $m$ objects in one of $a$ available colours. 
Since both sides are polynomials, this extends to negative integers as well.
For the case of an arbitrary binomial ring, invoke the Binomial Transfer Principle of \cite{BR}. 
\epr

\bth The Ariadne functor factors through the quotient category $\Laby_n$, producing a functor 
$$ A_n\colon \Laby_n\to \MSet_n.$$ \eth

\bpr It is clear that $A_n(P)=0$ when $|P|>n$. In order to prove that $A_n$ respects the relation 
$$ P=\sum_{\#A=P } \prod_{p\in P} \binom{\overline p}{\deg_A p} I_A,$$
it will be sufficient to establish that 
$$ A_n\left(Q\cup\spass{u}{v}{a}\right) 
= \sum_{k=1}^\infty \binom ak A_n\left( Q\cup \bigcup_k \spass{u}{v}{1} \right) $$
for any maze $Q$; one then performs induction on the passages of $P$.
This formula follows from a suitable application of the lemma, by which  
$$
A_n\left(Q\cup\spass{u}{v}{a}\right) = \sum_{m=1}^\infty 
\sum_{\substack{\#B=Q \\ |B|=n-m}} 
\left( \prod_{[q\colon x\to y]\in Q} \overline{q}^{\deg_B q}\col{x}{y}^{[\deg_B q]} \right)a^m\col{u}{v}
$$
is brought to co-incide with
\begin{align*}
& A_n\left(Q\cup\bigcup_k\spass{u}{v}{1}\right) \\ 
& \quad = \sum_{m=1}^\infty 
\sum_{\substack{\#B=Q \\ |B|=n-m}} \sum_{\substack{g_1+\cdots+g_k=m \\ g_i\in\Z^+}}
\left( \prod_{[q\colon x\to y]\in Q} \overline{q}^{\deg_B q}\col{x}{y}^{[\deg_B q]} \right) 
\col{u}{v}^{[g_1]}\cdots  \col{u}{v}^{[g_k]}. 
\end{align*} 
\epr

\bth The Ariadne functor factors through the quotient category $\Laby^n$, producing a functor 
$$ A_n\colon \Laby^n\to \MSet_n.$$ \eth

\bpr  
$$
A_n(a\boxdot P) 
= \sum_{\substack{\#A=P \\ |A|=n}} 
\left( \prod_{[p\colon x\to y]\in P} (a\overline{p})^{\deg_A p}\col{x}{y}^{[\deg_A p]} \right) 
= a^nA_n(P).
$$ 
\epr

The process of factorisation has now terminated, for passing 
to the category $\Laby^n$ has the effect of making the Ariadne functor faithful, which is
easily inferred from Theorem \ref{T: Generate} below.

\subsection{Pure Mazes}

The Ariadne Functor sheds light on the internal structure of the Labyrinth Categories.

\bdf A maze of which all passages carry the label $1$, is called a \textbf{pure} maze. \edf

\bth \label{T: Pure} The pure mazes are linearly independent in the category $\Laby$. \eth

\bpr 
Suppose we have a relation 
$$ \sum_j a_jP_j  = 0, $$
where the $P_j$ are distinct pure mazes in $\Laby(X,Y)$, for some $a_j\in \B$. Suppose all mazes have cardinality 
at least $n$. Applying the  $n$'th Ariadne functor will kill all mazes of cardinality greater than $n$, 
and the end result will be
$$ \sum_{|P_j|=n} a_jA_n(P_j) = 0. $$
Since the $P_j$ are distinct pure mazes, the $A_n(P_j)$ will all be distinct multations, which 
are linearly independent in $\MSet_n$. 
Hence all those $a_j=0$. The claim now follows by induction.
\epr

\bth \label{T: Basis} The pure mazes constitute a basis for the category $\Laby_n$. \eth

\bpr Linear independence goes through exactly as before. 
From the defining equations for $\Laby_n$, it follows that any maze will reduce to pure ones. \epr

\bth \label{T: Generate} 
The pure mazes with exactly $n$ passages are linearly independent in the category 
$\Laby^n$ and generate the category over $\Q\otimes_\Z\B$. \eth

\bpr Linear independence goes through as before. The defining equation for $\Laby^n$ can be written 
$$a^n P = \binom{a}{P}P + \sum_{\substack{\#A=P \\ |P|<|A|\leq n}} \binom{a}{A} A.$$ 
Since $a^n\neq \binom {a}{P}$ if $|P|<n$, such a $P$ may be expressed in terms of mazes with more passages, provided division by integers be permissible. \epr 

\bth \label{T: BZ} The following categories are isomorphic: 
$$ {}_\B\Laby_n \cong \B\otimes_\Z {}_\Z\Laby_n \qquad\qquad {}_\B\Laby^n \cong \B\otimes_\Z {}_\Z\Laby^n $$ 
\eth

\bpr
The first equation is a quick corollary of Theorem \ref{T: Basis}. 
Let us prove the second. By Theorem \ref{T: Generate}, $\Laby^n$ is torsion-free. 
From the definition of $\Laby^n$,  
$$ \Q\otimes_\Z {}_\B\Laby^n  \cong \Q \otimes_\Z \B \otimes_\Z  {}_\Z\Laby^n.$$
The category ${}_\B\Laby^n$ (canonically) embeds in the former and $\B\otimes_\Z {}_\Z\Laby^n$ in the latter. 
\epr

Because the Ariadne functor embeds ${}_\Z\Laby^n$ in ${}_\Z\MSet_n$, which is free, it follows that ${}_\Z\Laby^n$ 
is free as well. By the isomorphism above, ${}_\B\Laby^n$ will then be free for arbitrary $\B$, 
though it seems to possess no preferred basis.

\subsection{The Theseus Functor}

It is the purpose of the present section to adjoin objects to the category $\Q\otimes_{\Z}\Laby^n$ to the effect 
that certain sets 
split into direct sums. This will enable us to define an inverse to the Ariadne functor. 

\blem \label{L: Zero polynomial} 
Let $M$ be a torsion-free module, and let $p(x)\in M\otimes \B[x]$. Then $p=0$ if and only if 
$p(a)=0$ for all integers $a$. \elem

\blem \label{L: Split} In the category $\Q\otimes_\Z \Laby^n$, the following equation holds, for any set $X$: 
$$ I_X = \sum_{\substack{\#S=I_X \\ |S|=n}} \frac{1}{\deg S} S. $$ \elem

\bpr Use Lemma \ref{L: Zero polynomial}. Identify the co-efficient of $a^n$ in the defining equation for $\Laby^n$: 
$$a^nI_X = \sum_{\substack{\#S=P \\ |S|\leq n}} \binom{a}{S} S.$$ \epr

The mazes $\frac{1}{\deg S} S$, occurring in the lemma, satisfy 
$$ \frac{1}{\deg S} S \circ \frac{1}{\deg T} T  = 
\bca 0 & \text{if $S\neq T$,} \\ \frac{1}{\deg S} S & \text{if $S=T$.} \eca $$
They can thus be said to form a direct sum system, although the \emph{objects} themselves 
of the system do not exist. 
There is a simple remedy for this:  
adjoin to the category $\Q\otimes_\Z \Laby^n$ an object 
$ \Im \frac{1}{\deg S} S $ for each such $S$. 
The category $\Q\otimes_\Z \Laby^n$ is not, however, the minimal localisation of $\Laby^n$ for which 
this procedure makes sense. 

Let us say that passages $p\colon x\to y$ and $q\colon x\to y$, sharing starting and ending points, are 
\textbf{parallel}, and that a \textbf{simple} maze is one which contains no (pairs of) parallel passages. 
By means of the Labyrinth Axioms I and II, any maze can be written as the sum of simple mazes. 

It may be verified that the mazes $\frac{1}{\deg_P A} A$, 
where $|A|=n$ and $\#A=P$ for some pure and simple maze $P$, form a basis for a subcategory $C$
of $\Q\otimes \Laby^n$ which contains $\Laby^n$. 

\bdf \label{D: Laby+} To the category $C$, adjoin an object 
$ \Im \frac{1}{\deg S} S$ for each maze $S$, which, as a multi-set, is supported in some maze $I_X$, 
and denote the resulting category by $\Laby^{\oplus n}$.\edf

By Lemma \ref{L: Split}, the set $X$ will now split into components: 
$$ X = \bigoplus_{\substack{\#S=I_X \\ |S|=n}} \Im \frac{1}{\deg S} S. $$

\bex \label{E: Laby^2} It is a consequence of the equation 
$$ \begin{bmatrix} \xymatrix{\ast \upar[r]^1 \downar[r]_1 & \ast } \end{bmatrix} 
= 2 \begin{bmatrix} \xymatrix{\ast \ar[r]^1 & \ast } \end{bmatrix} $$
that no localisation will be required in the case $n=2$, so that, anomalously,
the categories $\Laby^2$ and $\Laby^{\oplus 2}$ are isomorphic. 
\eex

\bex Consider
$$ 
\xymatrixrowsep{0.2pc} \xymatrixcolsep{1pc} 
P=  \begin{bmatrix}  \xymatrix{
1 \upar[r]^1 \downar[r]_1  & 1 \\
2 \ar[r]_1 & 2
} \end{bmatrix} \qquad\text{and}\qquad 
Q = \begin{bmatrix}  \xymatrix{
1 \ar[r]^1 & 1 \\
2 \upar[r]^1 \downar[r]_1 & 2 
} \end{bmatrix}.
$$
In $\Q\otimes_{\Z}\Laby^3$, the equations $2I_{\{1,2\}}=P+Q$ and
$PQ=QP=0$ hold; hence, the mazes $\frac12 P$ and $\frac12 Q$ would form a direct sum system, 
were it not for the non-existence of the desired objects.  
By adjoining these, we obtain the refined category $\Laby^{\oplus 3}$, 
in which the set $\{1,2\}$ splits up into two:
$$ \{1,2\} = \Im \frac12 P \oplus \Im \frac 12 Q. $$ 
\eex

\bdf The $n$'th \textbf{Theseus functor} $$T_n\colon \MSet_n\to\Laby^{\oplus n}$$ is given by 
the following formul\ae: 
\begin{gather*} 
A\mapsto \Im \frac{1}{\deg A} \bigcup_{a\in A} \spass{a}{a}{1} \\
\mu\mapsto \frac{1}{\deg\mu} \bigcup_{(a,b)\in\mu} \spass{a}{b}{1}.
\end{gather*} \edf

It should be clear that this is indeed a 
(linear) functor, as composition in both categories is effectuated by ``summing over all possibilities''. 

\bth \label{T: Theseus} There is an isomorphism of categories: 
$$ \xymatrix{\Laby^{\oplus n} \upar[r]^{A_n} & \upar[l]^{T_n}\MSet_n } $$ \eth

\bpr We first shew that the Ariadne functor actually factors through the category $\Laby^{\oplus n}$. 
Suppose $P$ is a pure and simple maze, and let $A$ be supported in $P$, with $|A|=n$. Then 
$$ A_n(A) = \prod_{p=[1\colon x\to y]\in A} \col{x}{y} = 
\deg A \cdot \prod_{p=[1\colon x\to y]\in P} \col{x}{y}^{[\deg_A p]},$$ 
and hence we may extend $A_n$ by letting
$$ \frac{1}{\deg_P A} A \mapsto \prod_{p=[1\colon x\to y]\in P} \col{x}{y}^{[\deg_A p]}.$$ 

Moreover, $A_n$ maps the ``virtual'' biproduct system 
$$ I_X = \sum_{\substack{\#B=X \\ |B|=n}} \frac{1}{\deg B} \bigcup_{b\in B} \spass{b}{b}{1} $$
in $\Laby^n$ onto the ``real'' biproduct system 
$$ \sum_{\substack{\#B=X \\ |B|=n}} \iota_B$$
in $\MSet_n$, and we may consequently extend $A_n$ to $\Laby^{\oplus n}$ by defining
$$ A_n\left(\Im \frac{1}{\deg B} \bigcup_{b\in B} \spass{b}{b}{1} \right) = B.$$ 

It is now easy to see that $T_n$ and $A_n$ are inverse to each other. 
\epr

\subsection{The Category of Surjections} 
For reference, we devote this paragraph to investigating the connexion between the Labyrinth Category 
and the Category of Surjections  
explored in \cite{BDFP}. 
For want of space, we merely sketch the relevant constructions. 
 The reader anxious to learn the full details is referred to Section 5 of \emph{loc. cit.}

Let $C$ be a category possessing weak pull-backs. 
For two objects $X,Y\in C$, 
a \textbf{correspondence} or \textbf{span}
 from $X$ to $Y$ is a diagram 
$$ \xymatrixcolsep{1pc} \xymatrix{Y & \ar[l] U \ar[r] & X} $$
in $C$, intended to be read from right to left. 
Construct an additive category $\widehat C$ in the following way. Its objects are those of $C$. Its arrows 
are formal sums of correspondences of $C$ (identified under an obvious equivalence relation), living
in the free abelian group they generate. 
The composition of two correspondences is found by summing over weak pull-backs. 
The category $\widehat{C}$ is called the \textbf{category of correspondences} on $C$. 

It will now be observed that the category $\Sur$ of finite sets and \emph{surjections} possesses weak pull-backs. 
It is thus possible to build the category $\NQO$ of surjection correspondences. We form a quotient 
category $\NQO_n$ by forcing   
$$ \left[ \xymatrixcolsep{1pc} \xymatrix{Y & \ar[l] U \ar[r] & X} \right] = 0\quad\text{whenever $|U|>n$.} $$ 

\begin{theorem*} \label{T: BDFP} $$ \NQO_n \cong {}_\Z\Laby_n. $$ \end{theorem*}

\bpr Define the functor $$\Xi\colon \NQO_n\to {}_\Z\Laby_n$$ to be the identity on objects. 
The correspondence 
$$ \xymatrixcolsep{2pc} \phi= \left[ \xymatrix{Y & \ar[l]_-{\phi_\ast} U \ar[r]^-{\phi^\ast} & X} \right] $$
in $\NQO_n$ shall map to the pure maze $X\to Y$ with exactly
$$ \left|(\phi^\ast,\phi_\ast)^{-1}(x,y)\right| $$ passages going $X\ni x\to y\in Y$.

A simple and straight-forward calculation confirms that this yields a functor, which is inversible
because the pure mazes form a basis. 
\epr

Curiously enough,  
the categories $\NQO$ and ${}_\Z\Laby$ are themselves not isomorphic.
This stems from the fact that ${}_\Z\Laby$ encodes functors from the category of free abelian groups, 
while $\NQO$ was built to encode functors from 
the category of free commutative \emph{monoids}. These functor categories are not originally equivalent, 
but they will be, once polynomiality is brought into effect.

\section{The Combinatorics of Functors}

\subsection{Module Functors}
Our first aim is to shew how general module functors are encoded by the Labyrinth Category. 
The base ring $\B$ is presumed unital only, not necessarily commutative. 
It will be expedient to point out that 
finitely generated, free modules are automatically bimodules.
For homomorphisms, however, the left--right distinction is essential, and we hereby declare all 
maps under consideration to be 
\emph{right module homomorphisms} (and hence left multiplication by a matrix of ring elements). 

Let $$ \sigma_{yx}\colon\B^{X}\to\B^{Y}, \quad x\in X,\ y\in Y$$ denote the canonical \emph{transportations}, 
mapping a $1$ in position $x$ to a $1$ in position $y$, and everything else to $0$.

Moreover, let $\ce_X F(\B)$ denote the cross-effect of rank $X$ of the functor $F$, evaluated on $|X|$ copies 
of the module $\B$. Similarly, when $\zeta\colon F\to G$ is a natural transformation, let 
$$(\ce_X \zeta)_\B\colon \ce_X F(\B)\to\ce_X G(\B)$$ denote the evaluation of $\ce_X\zeta$ on $|X|$ copies of $\B$.

We propose a study of the functor $$\Phi\colon \MFun \to \Lin(\Laby,\Mod),$$ defined as follows. 
Given a module functor $F\colon \XMod\to \Mod$, the corresponding labyrinth functor 
$\Phi(F)\colon\Laby\to\Mod$ should take:
\begin{gather*}
X \mapsto \ce_X F(\B) \\
[P\colon X\to Y] \mapsto \left[ F\left(\De_{[p\colon x\to y]\in P} \overline{p}\sigma_{yx}\right) \colon \ce_X F(\B) \to \ce_Y F(\B)  \right]. 
\end{gather*}
Of course, one ought to restrict the action to the appropriate cross-effects, but this turns out to be an 
unnecessary caution:

\blem The map 
$$\upsilon=F\left(\De_{[p\colon x\to y]\in P} \overline{p}\sigma_{yx}\right) \colon F(\B^X) \to F(\B^Y)$$ 
is in fact a map $\ce_X F(\B) \to \ce_Y F(\B)$, in the sense that there is a commutative diagram: 
$$  \xymatrix{ 
F(\B^X) \ar[d] \ar[r]^{\upsilon} & F(\B^Y) \\
\ce_X F(\B) \ar[r]_{\upsilon} & \ce_Y F(\B) \ar[u]
}$$ \elem 

\blem $\Phi(F)$ is a functor $\Laby\to\Mod$. \elem

\bpr That $\Phi(F)$ respects the relations in $\Laby$ follows from
$$\Phi(F)\left( P\cup\{ \xymatrix{x\ar[r]^0 & y }\} \right)= F(\cdots \de 0) = 0$$
and
\begin{align*} 
& \xymatrixcolsep{1.5pc}
\Phi(F)\left(P\cup\{ \xymatrix{x\ar[r]^{a+b} & y }\} \right) = F(\cdots \de (a+b)\sigma_{yx}) \\
&\quad = F(\cdots \de a\sigma_{yx}) + F(\cdots \de b\sigma_{yx}) + F(\cdots \de a\sigma_{yx}\de b\sigma_{yx}) \\
&\quad \xymatrixcolsep{0.95pc}
 =\Phi(F)\left( P\cup\{ \xymatrix{x\ar[r]^a & y }\} \right) 
+ \Phi(F)\left(  P\cup\{ \xymatrix{x\ar[r]^b & y }\} \right)  
+ \Phi(F)\left( P\cup\{ \xymatrix{x\upar[r]^a \downar[r]_b & y }\} \right).
\end{align*}
Functoriality of $\Phi(F)$ is a consequence of the Deviation Formula and the definition of maze composition. \epr

Let $\zeta\colon F\to G$ be a natural transformation. 
Define $\Phi(\zeta)\colon \Phi(F)\to\Phi(G)$ by restriction to the appropriate cross-effects: 
$$ \Phi(\zeta)_X= (\ce_X \zeta)_\B \colon \ce_X F(\B) \to \ce_X G(\B).$$

\blem $\Phi$ is a functor $$\MFun \to \Lin(\Laby,\Mod).$$ \elem

\bpr Follows from the functoriality of $\ce_X$. \epr

\blem $\Phi$ is fully faithful. \elem 

\bpr The natural transformation $\zeta$ can be uniquely re-assembled from its components $\ce_X\zeta$. \epr

We now construct the inverse of $\Phi$.  Given a labyrinth functor $H\colon \Laby\to \Mod,$ 
the corresponding module functor $\Phi^{-1}(H)$ should take:  
\begin{gather*} 
\B^A \mapsto \bigoplus_{X\subseteq A} H(X) \\
\left[\sum_{\substack{a\in A \\ b\in B}} s_{ba}\sigma_{ba}\colon \B^A\to \B^B \right]
\mapsto \sum_{K\subseteq B\times A} H\left( \bigcup_{(b,a)\in K} \Set{ \ssap{b}{a}{s_{ba}} } \right).
\end{gather*}

\blem $\Phi^{-1}(H)$ is a module functor. \elem

\bpr Functoriality is established thus:
\begin{align*} 
& \Phi^{-1}(H)\left(\sum_{\substack{b\in B \\ c\in C}} s_{cb}\sigma_{cb}\right) \circ 
\Phi^{-1}(H)\left(\sum_{\substack{a\in A \\ b\in B}} t_{ba}\sigma_{ba}\right) \\
&= \sum_{I\subseteq C\times B} H\left( \bigcup_{(c,b)\in I} \Set{ \ssap{c}{b}{s_{cb}} } \right) \circ
\sum_{J\subseteq B\times A} H\left( \bigcup_{(b,a)\in J} \Set{ \ssap{b}{a}{t_{ba}} } \right) \\ 
&=  \sum_{I\subseteq C\times B} \sum_{J\subseteq B\times A}  
\sum_{K\sqsubseteq \{((c,b),(b,a))\in I\times J\} } 
H\left( \bigcup_{((c,b),(b,a))\in K} \Set{ \ssap{c}{a}{s_{cb}t_{ba}} } \right) \\ 
&= \sum_{M\subseteq C\times B\times A} 
H\left( \bigcup_{(c,b,a)\in M} \Set{ \ssap{c}{a}{s_{cb}t_{ba}} } \right) \\ 
&= \sum_{L\subseteq C\times A} H\left( \bigcup_{(c,a)\in L} 
\Set{ \xymatrixcolsep{3pc}  \ssap{c}{a}{\sum_{b\in B} s_{cb}t_{ba}} } \right) 
= \Phi^{-1}(H)\left(\sum_{\substack{a\in A \\ c\in C}} \left(\sum_{b\in B} s_{cb}t_{ba}\right)\sigma_{ca}\right),
\end{align*}
where Theorem \ref{T: Maze formulae} was used in the fifth step. \epr

\begin{lemma*} $$\Phi(\Phi^{-1}(H))=H.$$\end{lemma*}

\bpr Let $P\colon X\to Y$ be a maze. We calculate:  
\begin{align*} 
\Phi(\Phi^{-1}(H))(P) &= \Phi^{-1}(H)\left(\De_{[p\colon x\to y]\in P} \overline{p}\sigma_{yx}\right) 
= \sum_{S\subseteq P} (-1)^{|P|-|S|} \Phi^{-1}(H)\left(\sum_{p\in S} \overline{p} \sigma_{yx} \right) \\
&= \sum_{S\subseteq P} (-1)^{|P|-|S|} \sum_{K\subseteq S} H(K) 
= \sum_{K\subseteq P} H(K)  \sum_{K\subseteq S\subseteq P} (-1)^{|P|-|S|} = H(P).
\end{align*}
 \epr

Assembling these results, we conclude that $\Phi$ is fully faithful and representative, and thence 
we obtain our main theorem. 

\bth[Labyrinth of Fun]The functor \label{T: LoF}
$$ \Phi\colon \MFun\to \Lin(\Laby,\Mod), $$
where $\Phi(F)\colon \Laby\to\Mod$ takes
\begin{gather*}
X \mapsto \ce_X F(\B) \\
[P\colon X\to Y] \mapsto \left[ F\left(\De_{[p\colon x\to y]\in P} \overline{p}\sigma_{yx}\right) 
\colon \ce_X F(\B) \to \ce_Y F(\B) \right] , 
\end{gather*}
is an equivalence of categories. \eth

\subsection{Polynomial Functors}

Since mazes correspond to deviations, this very simple characterisation of polynomiality should come as no surprise.  

\bth \label{T: LoF Pol} The module functor $F$ is polynomial of degree $n$ if and only if 
$\Phi(F)$ vanishes on sets with more than $n$ elements; or, equivalently, on mazes with more than 
$n$ passages. \eth

\bpr Assume first that $F$ is polynomial of degree $n$. 
Since mazes with $n+1$ passages correspond to $n$'th deviations, 
$\Phi(F)$ will certainly vanish on mazes with more than $n$ passages.

Suppose now, conversely, that $\Phi(F)$ annihilates mazes with more than $n$ passages. 
Then each cross-effect of $F$ of rank exceeding $n$ will vanish, and $F$ is polynomial of degree $n$. 
\epr

\subsection{Numerical Functors}

We now investigate how to interpret numericality in the labyrinthine setting. 
The base ring $\B$ will of course be assumed binomial.

\blem Let $r\in\B$, and let $w_1,\dots,w_q$ be natural numbers. Then 
$$ \prod_{j=1}^q \binom{r}{w_j} =  
\sum_{m=0}^\infty \binom{r}{m} \sum_{k=0}^{m} (-1)^{m-k} \binom{m}{k} \prod_{j=1}^q \binom{k}{w_j} . $$ \elem

\bpr 
When $r$ is an integer, both sides of the given equality count the number of ways to choose subsets 
$W_1,\dots,W_q\subseteq [r]$ with $|W_i|=w_i$. This is because, for a fixed subset $S\subseteq [r]$ with $|S|=m$,
there are, by the Principle of Inclusion and Exclusion, exactly 
$$\sum_{k=0}^{m} (-1)^{m-k} \binom{m}{k} \prod_{j=1}^q \binom{k}{w_j}$$ subsets 
$W_1,\dots,W_q\subseteq [r]$ such that $|W_i|=w_i$ and $\bigcup W_i=S$. 

When $r$ is an element of an arbitrary binomial ring, we invoke the Binomial Transfer Principle of \cite{BR}. \epr

Recall that, when $A$ is a maze (hence a multi-set), the symbol $I_A$ betokens the result of substituting $1$ 
for the labels of all the passages of $A$.

\bth \label{T: LoF Num} 
The module functor $F$ is numerical of degree $n$ if and only if $\Phi(F)$ vanishes on sets (or mazes) 
with more than $n$ elements, and, withal, satisfies the equation
$$ \Phi(F)(P) = \sum_{\#A=P } \prod_{p\in P} \binom{\overline p}{\deg_A p} \Phi(F)(I_A) $$
for any maze $P$; so that $\Phi(F)$ factors through $\Laby_n$. 
The functor $\Phi$ induces an equivalence of categories $$ \Num_n \to \Lin(\Laby_n,\Mod). $$ \eth

\bpr 
By the previous theorem, $\Phi(F)$ vanishing on sets with more than $n$ elements is equivalent to polynomiality. 
It is then clear from  
Theorem 10 of \cite{PF} that numerical functors satisfy the stated requirements.  

Suppose now, conversely, that $\Phi(F)$ satisfies the conditions of the theorem. 
We wish to use Criterion A$'$ in Theorem 9 of 
\cite{PF} and thus seek to evaluate 
\begin{align*}
F(r\cdot 1_{\B^n}) &= \sum_{P\subseteq r\boxdot I_{[n]}} \Phi(F)(P) 
= \sum_{J\subseteq [n]}\Phi(F)\left( \bigcup_{j\in J} \spass{j}{j}{r} \right) 
\end{align*}
and
\begin{align*}
\sum_{m=0}^\infty \binom rm F\left(\De_m 1_{\B^n}\right) 
&= \sum_{m=0}^\infty \binom rm \sum_{k=0}^m (-1)^{m-k} \binom mk F( k\cdot 1_{\B^n} ) \\
&= \xymatrixcolsep{1.55pc}
\sum_{m=0}^\infty \binom rm \sum_{k=0}^m (-1)^{m-k} \binom mk \sum_{J\subseteq [n]} 
\Phi(F)\left( \bigcup_{j\in J} \spass{j}{j}{k} \right).
\end{align*}
To verify the equality 
$$ F(r\cdot 1_{\B^n}) = \sum_{m=0}^\infty \binom rm F\left(\De_m 1_{\B^n}\right), $$
it will then be sufficient to check that the terms corresponding to a certain fixed $J$ are equal. 
There will be no loss of generality in considering the special case $J=[q]$ only. Our object will thus be to verify 
$$ \xymatrixcolsep{1.5pc}
\Phi(F)\left( \bigcup_{j\in [q]} \spass{j}{j}{r} \right) = 
\sum_{m=0}^\infty \binom rm \sum_{k=0}^m (-1)^{m-k} \binom mk \Phi(F)\left( \bigcup_{j\in [q]} \spass{j}{j}{k} \right).
$$
By assumption, the left-hand side equals 
$$ \sum_{\substack{w_1+\cdots+w_q\leq n \\ w_j\geq 1}} \prod_{j=1}^q \binom{r}{w_j} 
\Phi(F)\left(\bigcup_{j\in[q]}\bigcup_{w_j} \spass{j}{j}{1}\right) $$
and the right-hand side 
$$ \sum_{m=0}^\infty \binom rm \sum_{k=0}^m (-1)^{m-k} \binom mk \sum_{\substack{w_1+\cdots+w_q\leq n \\ w_j\geq 1}} 
\prod_{j=1}^q \binom{k}{w_j} \Phi(F)\left(\bigcup_{j\in [q]}\bigcup_{w_j} \spass{j}{j}{1} \right),$$
so that the equality follows after deployment of the lemma. \epr

\bex As a simple example, a labyrinth module $H$ corresponding to a cubical functor will 
satisfy the equation
\begin{align*}  
H\begin{bmatrix}  \xymatrix{ \ast \ar@/^1pc/[r]^a \ar@/^-1pc/[r]_b  & \ast  }  \end{bmatrix}
& = 
\binom a1 \binom b1  
H\begin{bmatrix}  \xymatrix{  \ast \ar@/^1pc/[r]^1 \ar@/^-1pc/[r]_1 & \ast }\end{bmatrix} \\
& + \binom a2 \binom b1  
H\begin{bmatrix}  \xymatrix{  \ast \ar@/^1pc/[r]^1 \ar@/^.3pc/[r]_1 \ar@/^-1pc/[r]_1 & \ast }\end{bmatrix}  
+ \binom a1 \binom b2 
H\begin{bmatrix}  \xymatrix{  \ast \ar@/^1pc/[r]^1 \ar@/^-.3pc/[r]^1 \ar@/^-1pc/[r]_1  & \ast }\end{bmatrix}.
\end{align*}
\eex

The main result of Baues, Dreckmann, Franjou, and Pirashvili in \cite{BDFP}: 
$${}_\Z\Num_n \cong \Lin(\NQO_n,{}_\Z\Mod),$$ 
will now be obtained as a simple corollary from the category equivalences exhibited in 
Theorems \ref{T: BDFP} and \ref{T: LoF Num}.

\subsection{Quasi-Homogeneous Functors}

\bth \label{T: LoF QHom} 
The module functor $F$ is quasi-homogeneous of degree $n$ if and only if $\Phi(F)$ satisfies the equation 
$$ \Phi(F)(a\boxdot P) = a^n \Phi(F)(P) $$
for any maze $P$ and scalar $a$; so that $\Phi(F)$ factors through $\Laby^n$.
The functor $\Phi$ induces an equivalence of categories $$ \QPol_n \to \Lin(\Laby^n,\Mod). $$ \eth

\bpr Let $F$ be quasi-homogeneous, and let $a\in\Q\otimes_\Z \B$. For any deviation, we have 
$$ F(a\alpha_1\de\cdots\de a\alpha_k) = a^n F(\alpha_1\de\cdots\de\alpha_k), $$
and we may calculate for a pure maze $P$: 
$$
a^n \Phi(F)(P) = a^n F\left(\De_{[p\colon x\to y]\in P} \overline{p}\sigma_{yx}\right) = 
F\left(\De_{[p\colon x\to y]\in P} a\overline{p}\sigma_{yx}\right) = \Phi(F)(a\boxdot P) .
$$
Conversely, assume $\Phi(F)$ factors via $\Laby^n$. Then, for any $k\in\N$, 
$$ a^nF(1_{\B^k}) = a^n\sum_{K\subseteq [k]} \Phi(F)(I_K) 
= \sum_{K\subseteq [k]} \Phi(F)(a\boxdot I_K) = F(a\cdot 1_{\B^k}), $$
 and $F$ is quasi-homogeneous. \epr

\subsection{Strict Polynomial Functors}
The following theorem was first formulated by Salomonsson in terms of Mackey functors, and later reformulated in 
\cite{X} using multations. 

\begin{reftheorem}[\cite{Pelle}, Theorem I.2.3; \cite{X}, Theorem 10.9] \label{T: LoF Hom} The functor 
$$ \Psi\colon \HPol_n \to \Lin(\MSet_n,\Mod), $$
where $$\Psi(F)\colon \MSet_n\to\Mod$$ takes
\begin{gather*}
A \mapsto \Im \pcom{F}{\pi}{A} \\
\left[\mu\colon A\to B\right] \mapsto \left[\pcom{F}{\sigma}{\mu}\colon \ce_A F(\B)\to \ce_B F(\B)\right],
\end{gather*}
is an equivalence of categories. \end{reftheorem}

It will be of interest to spell out a formula for the inverse. 
When $J\colon \allowbreak\MSet_n\to\Mod$, the functor $\Phi^{-1}(J)\colon\XMod\to\Mod$ is defined by 
\begin{gather*} 
\B^X \mapsto \bigoplus_{\substack{\#A\subseteq X \\ |A|=n}} J(A) \\
\left[\sum_{\substack{x\in X \\ y\in Y}} s_{yx}\sigma_{yx} \colon \B^X\to \B^Y\right] \mapsto 
\sum_{\substack{\#A\subseteq X,\ \#B\subseteq Y \\ |A|=|B|=n}} 
\sum_{\mu\colon A\to B} s^\mu J(\mu)
\end{gather*}
(where, of course, $X$ and $Y$ are sets, but $A$ and $B$ range over multi-sets).

\section{Numerical versus Strict Polynomial Functors} 
In this final section, we provide a comparison of the two strains of functors: numerical (polynomial) 
and strict polynomial.

\subsection{Quadratic Functors} We first propose to examine quadratic functors in detail, and thus seek to
fathom the structure of $\Num_2$. The key point lies in unravelling 
the structure of the category $\Laby_2$. It contains three non-isomorphic objects: $[0]$, $[1]$, and $[2]$. 
By Theorem \ref{T: Pure}, the pure mazes form a basis. Those are identity mazes, along with the  
 quadruple:  
$$ \xymatrixrowsep{0.2pc} \xymatrixcolsep{1pc}   
A=\begin{bmatrix}  \xymatrix{
1 \ar[r]^(0.7){1} \ar[dr]_(0.7){1}  & 1  \\
& 2 } \end{bmatrix}, \qquad
B = \begin{bmatrix} \xymatrix{
1 \ar[r]^(0.7){1} & 1   \\
2 \ar[ur]_(0.7){1} 
} \end{bmatrix}, \qquad
\xymatrixrowsep{0.2pc} \xymatrixcolsep{1pc} 
C = \begin{bmatrix} \xymatrix{
1 \upar[r]^{1} \downar[r]_{1} & 1  
} \end{bmatrix}, \qquad 
S=\begin{bmatrix} \xymatrix{
1 \ar[dr]_(0.8){1} & 1 \\
2 \ar[ur]^(0.8){1} & 2
} \end{bmatrix} .
$$
The (skeletal) structure of the category $\Laby_2$ is thus reduced to the following, promptly suggesting the nick-name \emph{dogegory}: 
$$ \xymatrix{
[0]\ar@(dl,ul)^{I} & [1] \curvear[r]^A \ar@(d,l)^{I} \ar@(l,u)^{C} & [2] \curvear[l]^B \ar@(u,r)^{S} \ar@(r,d)^{I}
} $$
An inspection of the multiplication table, given in Table \ref{Ta: Laby_2}, 
\begin{table} \centering
\begin{minipage}[b]{.45\textwidth}
  \centering
  \begin{tabular}{|c|cccc|}
\hline 
$\circ$ & $A$  & $B$   & $C$  & $S$  \\
\hline                        
$A$     & --   & $I+S$ & $2A$ & --   \\               
$B$     & $C$  & --    & --   & $B$  \\                 
$C$     & --   & $2B$  & $2C$ &  --  \\ 
$S$     & $A$  & --    & --   & $I$  \\ 
\hline 
  \end{tabular}
  \caption{Multiplication table for $\Laby_2$. \label{Ta: Laby_2}}
\end{minipage}\qquad
\begin{minipage}[b]{.45\textwidth}
  \centering
  \begin{tabular}{|c|ccc|}
\hline 
$\circ$      & $\alpha$  & $\beta$        & $\sigma$  \\
\hline
$\alpha$     & --        & $\iota+\sigma$ & --       \\               
$\beta$      & $2\iota$  & --             & $\beta$  \\                 
$\sigma$     & $\alpha$  & --             & $\iota$  \\ 
\hline 
  \end{tabular}
  \caption{Multiplication table for $\MSet_2$. \label{Ta: MSet_2}}
\end{minipage}
\end{table}
reveals that the mazes $A$, $B$, $C$, and $S$ are not algebraically independent, for 
$C=BA$ and $S=AB-I$. We thus recuperate the now classical classification 
of quadratic (integral) functors from \cite{Baues}: 

\bth A quadratic numerical functor is equivalent to a diagram of modules and homomorphisms as indicated, 
subject to the two relations:
\begin{gather*}
\beta\alpha\beta = 2\beta, \qquad \alpha\beta\alpha=2\alpha . 
\end{gather*}
$$ \xymatrix{
K & X \upar[r]^\alpha   & Y \upar[l]^\beta  
} $$ \label{Th: Num Bridge} \eth

To determine $\HPol_2$, we proceed similarly. The (skeletal) structure of the category $\MSet_2$ is:
$$ 
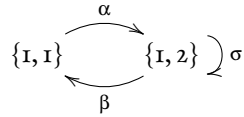
\xymatrix{
\{1,1\} \curvear[r]^\alpha & \{1,2\} \curvear[l]^\beta \ar@(ur,dr)^{\sigma} 
} $$
Every multation reduces to a linear combination of identity multations and the subsequent triplet, 
with multiplication given in Table \ref{Ta: MSet_2}:  
\begin{gather*} 
\alpha= \begin{bmatrix}
1 & 1 \\
1 & 2
\end{bmatrix} \qquad
\beta= \begin{bmatrix}
1 & 2 \\
1 & 1
\end{bmatrix} \qquad 
\sigma= \begin{bmatrix}
1 & 2 \\
2 & 1
\end{bmatrix}
\end{gather*} 

\bth 
A quadratic homogeneous functor is equivalent to a diagram of modules and homomorphisms as indicated, 
subject to the single relation:
\begin{gather*}
\beta\alpha=2.
\end{gather*}
$$ \xymatrix{X \upar[r]^\alpha   & Y \upar[l]^\beta  
} $$  \eth

We then obtain the following characterisation of homogeneous quadratic functors. 

\bth \label{T: Quadratic} Consider the labyrinthine description of a quadratic functor $F$: 
$$ \xymatrix{ K & X \upar[r]^\alpha & Y \upar[l]^\beta  } $$ 
The following conditions are equivalent: 
\balph
\item $F$ is quasi-homogeneous of degree $2$.
\item $F$ may be (uniquely) extended to a homogeneous quadratic functor.  
\item $K=0$ and $\beta\alpha = 2$.
\ealph\hfill \eth

\bpr
The equivalence of B and C follows from the two preceding theorems; that of A and B is a consequence of 
the isomorphisms 
$$ \MSet_2 \cong \Laby^{\oplus 2} \cong \Laby^2, $$
exhibited in Example \ref{E: Laby^2} and Theorem \ref{T: Theseus}. 
\epr

\bex \label{E: Frobenius} A wee example will serve to illustrate the theorem, and also to point out its subtlety. 
This same delicacy was observed in \cite{PF}, Example 1. 
Consider, over the ring $\Z$, the following labyrinth module: 
$$ \xymatrix{
0 & \Z/2\Z \upar[r]^-0 & 0 \upar[l]^-0  
} $$ 
Since it satisfies the conditions of the theorem, it has a unique structure of homogeneous quadratic functor 
--- viz.~the functor $G$ given by: 
$$ \xymatrix{
 \Z/2\Z \upar[r]^-0 & 0 \upar[l]^-0  
} $$ 
Yet it is possible to exhibit another strict polynomial structure $F$ on this same underlying functor, 
this one linear. As a multi-set module $$\bigoplus_{n=0}^2 \MSet_n \to \Mod,$$ it is given by the formul\ae
\begin{align*}
\{\} & \mapsto 0 		& \{1,1\} &\mapsto 0 \\
\{1\}& \mapsto \Z/2		& \{1,2\} &\mapsto 0.
\end{align*}
On the level of functors, $F=\Z/2\Z\otimes-$ and
$G=  (\Z/2\Z)^{(1)} \otimes - $, the symbol ${}^{(1)}$ indicating Frobenius twist. 
\eex

\subsection{The Ariadne Thread}
By the theory previously wrought out (Theorems \ref{T: LoF Num} and \ref{T: LoF Hom}),  
there are category equivalences 
$$ \Phi\colon\Num_n\to\Lin(\Laby_n,\Mod), \qquad \Psi\colon\HPol_n\to\Lin(\MSet_n,\Mod). $$
Connecting the combinatorial categories is the Ariadne functor $A_n\colon\Laby_n\to\MSet_n$.

\begin{theorem}[The Ariadne Thread] \label{T: Ariadne} Pre-composition with the Ariadne functor $A_n$ begets 
the forgetful functor $$\HPol_n\to\Num_n,$$ so that $$ \Phi \circ \Psi^{-1} = (A_n)^\ast. $$ \end{theorem}

\bpr Let $J\colon \MSet_n\to \Mod$ 
and let $P\colon X\to Y$ be a maze. We compute:
\begin{align*} 
\Phi\Psi^{-1}(J)(P) &= \Psi^{-1}(J)\left(\De_{[p\colon x\to y]\in P}\overline{p}\sigma_{yx}\right) 
= \sum_{I\subseteq P} (-1)^{|P|-|I|} \Psi^{-1}(J)\left(\sum_{p\in I} \overline{p}\sigma_{yx}\right) \\ 
&=\sum_{I\subseteq P} (-1)^{|P|-|I|} \sum_{\substack{\#A\subseteq I \\ |A|=n}} 
J\left(\prod_{p\in \#A} \overline{p}^{\deg_A p}\col{x}{y}^{[\deg_A p]}\right)  \\
&= \sum_{\substack{\#A\subseteq P \\ |A|=n}} J\left(\prod_{p\in \#A} \overline{p}^{\deg_A p}\col{x}{y}^{[\deg_A p]}\right)  
\sum_{\#A\subseteq I\subseteq P} (-1)^{|P|-|I|}  \\
&= \sum_{\substack{\#A= P \\ |A|=n}} J\left(\prod_{p\in \#A} \overline{p}^{\deg_A p}\col{x}{y}^{[\deg_A p]}\right)  = JA_n(P),
\end{align*}
whence $\Phi\Psi^{-1}(J) = J A_n$, as required.
\epr

\subsection{The Polynomial Functor Theorem}
Let us finally record the exact obstruction for a numerical functor, of arbitrary degree, to be strict polynomial. 

\bth[The Polynomial Functor Theorem] \label{T: PFT} Let 
$F$ be a quasi\hyp homogeneous functor of degree $n$, corresponding to the labyrinth module 
$H\colon\Laby^n\to\Mod$. 
 Imposing the structure of homogeneous functor upon $F$ is equivalent to  
 exhibiting a factorisation of $H$ through $\Laby^{\oplus n}$. 
\eth

\bpr It will suffice to recall the isomorphism $\Laby^{\oplus n}\cong\MSet_n$ from Theorem \ref{T: Theseus}. \epr

A word of admonishment: the existence of a factorisation of $H$ through $\Laby^{\oplus n}$ 
does not imply its uniqueness. As was seen in Example \ref{E: Frobenius}, there are, in general, 
many strict polynomial structures on the same functor, even of different degrees. 

\bex We wish to draw the reader's attention to one particular case, when factorisation always takes place. 
If $\B$ is a $\Q$-algebra, the categories
$$ \Laby^n=\Laby^{\oplus n} \cong \MSet_n$$
are isomorphic. This mirrors the well-known 
fact that, over a $\Q$-algebra, numerical and strict polynomial functors co-incide. 
\eex 

\bex For \emph{affine} functors (degree $0$ and $1$), there is no discrepancy between numerical and strict 
polynomial functors. This will no longer be the case in higher degrees. 
Yet, the quadratic case will retain some regularity, in that, according to Theorem~\ref{T: Quadratic}, 
any \emph{quasi-homogeneous} functor necessarily admits a unique homogeneous structure. 
\eex

\end{document}